\documentclass{report}

\newcommand{\comb}[2]{{#2 \choose #1}}

\usepackage[francais]{babel}
\usepackage[latin1]{inputenc}

\newcommand{\citer}[2]{Page #1, \S\,#2}
\newcommand{\citerp}[2]{\S\,#2}

\newcommand{\score}[2]{($#1$, $#2$)}

\newcommand{\w}[1]{\makebox[.4cm][c]{$#1$}}
\newcommand{\x}{\small}
\newcommand{\z}[1]
{{\makebox[.5cm][c]{\x\rule[-.3cm]{0cm}{0.9cm}#1}}}
\newcommand{\y}[2]{\mbox{\x$\begin{array}{c}#1\\#2\end{array}$}}

\title{Pascal, Fermat et la géométrie du hasard}
\author{Nicolas Trotignon \\ \vspace{.4cm}{\small
Mémoire réalisé à l'IUFM de Créteil}\\\vspace{.0cm}{\small
Sous la direction d'Evelyne Barbin}}  
\date{1er juin 1998}

\renewcommand{\thesection}{\arabic{section}}

\begin{document}

\maketitle
\thispagestyle{empty}
\markright{Table des matières}

\tableofcontents

\newpage
\section{Remeciements}
\setcounter{page}{1}

\vspace{5cm}

Je tiens à remercier Evelyne Barbin qui a dirigé ce travail
et a su me faire profiter de sa grande connaissance de l'histoire des
mathématiques. Je remercie également Christelle Petit d'avoir
patiemment relu mon texte. 

\vspace{1cm}

Vincent Gérard et Pierre-Étienne Moreau m'ont aidé à préciser
quelques indications bibliographiques.

\vspace{1cm}

La mise en forme de ce document a été réalisée sans
l'aide d'aucun logiciel commercial.

\newpage
\section{Introduction}

\subsection*{Le problème des partis}

``Deux camps jouent à la balle;~chaque manche est de dix points et il
faut 60 points pour gagner le jeu;~la mise totale est de dix
ducats. Il arrive que, pour quelque raison accidentelle, le jeu ne
puisse s'achever. On demande ce que touche de la mise totale chacun
des deux camps, lorsque l'un a 50 points et l'autre 20
points''%
 \footnote{
Luca Pacioli, \emph{Suma de arithmetica} (1509), d'après un
document polycopié de l'IREM de Paris~VII.
}.

L'énoncé ci-dessus, dû à Pacioli,  est la première version de ce
que 145 années 
plus tard Pascal appellera le problème des partis
\footnote{ Je signale  que dans les textes de  Pascal, il est question
de ``parties''  et de ``partis''.  ``Partie'' doit  s'entendre dans le
sens courant de ``durée d'un  jeu, à l'issue de laquelle sont désignés
gagnants et perdants''. ``Parti'', nom masculin vieilli, signifie ``ce
qu'une  personne  a  pour  sa  part'',  au moment  où  la  partie  est
interrompue  dans  le  cas  qui  nous occupe.   Ces  définitions  sont
empruntées au \emph{Petit Robert}.}.
Il réapparaît tout au long du XVI\ieme\ siècle sous la plume de
différents 
auteurs italiens, sans qu'aucun ne parvienne à en donner une solution
satisfaisante. Ernest Coumet évoque notamment les contributions de
Pacioli (1494), de Tartaglia (1578), de Forestani (1603), de Cardan
(1539) 
et de Peverone (1558)%
\footnote{
Ernest Coumet, \emph{Le problème des partis avant Pascal
\emph{in} Archives internationales d'histoire des sciences},
numéro 18, 1965, fascicule 72--73, pages 245\,--\,272.}.
En 1654, le
problème est enfin résolu par Pascal, puis,
par Fermat, ce qui marque la naissance de la théorie
des probabilités%
\footnote{
Le mot `probabilité' est ici anachronique. Il a une histoire
touffue. Je mentionnerai simplement qu'au XVII\ieme\ siècle, il
désignait ``la chance qu'une chose a d'être vraie'' et, en théologie,
une ``doctrine fondée sur les opinions 
probables''\footnote{Voir le dictionaire historique de la langue
française},
de laquelle, par une ironie de l'histoire, Pascal fut
l'adversaire.

Dans son acception mathématique, le mot probabilité n'a commencé à être
utilisé qu'au début du XVIII\ieme\ siècle.

Je constate que l'ordinateur refuse d'envisager la possibilité d'une
note en bas de page référencée dans une autre note (l'ordinateur n'est
pas borné : \emph{c'est} une borne). La note 5 renvoie au
\emph{Dictionnaire historique de la langue française} chez Robert.
}.

Pascal avait la conviction qu'il inaugurait là une nouvelle ère de
l'histoire des 
mathématiques et de la pensée :
dans son \emph{Adresse à l'Académie 
Parisienne}, en 1654, il donne un aperçu de ses travaux
achevés ou en 
projet, et annonce ``un traité tout à
fait nouveau, d'une matière absolument inexplorée jusqu'ici'', qui
peut ``s'arroger à bon droit ce titre stupéfiant : \emph{la géométrie
du hasard}%
 \footnote{Le traité s'appellera finalement \emph{Traité du
 triangle arithmétique}.}''.
``Stupéfiant'' car auparavant le hasard
semblait par essence irréductible à toute considération rationnelle~---%
~irréductibilité qui incitait au doute quant à la possibilité même de
résoudre le problème des partis. Ainsi Tartaglia pouvait-il écrire :  
``La résolution d'une telle question est davantage d'ordre judiciaire
que rationnel, et de quelque manière qu'on veuille la résoudre, on y
trouvera sujet à litige''%
 \footnote{Citation empruntée à E. Coumet,
 \emph{op. cit.}}.
(Cette situation demeure assez rare en histoire des mathématiques pour
mériter une parenthèse, les mathématiciens
ayant  d'ordinaire une confiance plutôt
excessive en la portée de concepts  qui montrent leurs
limites après coup~:~nombres rationnels, constructions à la règle
et au compas, résolution par radicaux des équations, géométrie
Euclidienne, systèmes formels \ldots). 

On comprend finalement que la résolution du problème des partis, plus qu'un
simple progrès technique des mathématiques, est une véritable révolution dans
notre conception du hasard. Et si le problème des
partis a une histoire avant Pascal\footnote{Voir à ce sujet Ernest
Coumet, \emph{op. cit.}},
il reste que le basculement définitif de l'erreur dans la vérité, opéré en
1654, marque une rupture nette et rapide dans le progrès des idées.
En dépit des prédécesseurs de Pascal, l'invention 
de la géométrie du hasard est un événement bref dont témoigne un
corpus de textes bien délimité, qu'il nous reste à examiner brièvement
avant de poursuivre.

\vspace{.4cm}
\subsection*{Les sources utilisées}

En 1654, le problème des partis resurgit dans des conditions
assez obscures. Le chevalier de Méré le propose à Pascal\footnote{``M. le
chevalier de Méré, qui est celui qui m'a proposé ces questions'',
lettre de Pascal à Fermat du 29 juillet 1654.}, qui le résout et
le soumet à son tour (par l'entremise de M. de Carcavi) à Fermat, qui
le résout lui aussi. Les textes
qui nous sont parvenus sont tous postérieurs à ces échanges (à l'exception
peut-être d'une lettre non datée). Ils consistent en une
correspondance entre Pascal, Fermat, Huygens et M. de Carcavi, ainsi qu'en
quelques traités donc voici une brève description à laquelle on pourra
se référer pour situer tel ou tel extrait donné par la
suite%
\footnote{Pour plus de détails, voir en bibliographie.} 
(un encart détachable reprend l'essentiel de cette liste)
:

\begin{itemize}
\item Correspondance :

  \begin{itemize}
	\item 
          {\bf Une lettre non datée de Fermat à Pascal.}
	
	Assez obscure, cette lettre concerne le jeu de
	dés. Nous l'appellerons `la lettre
	non datée'.

	 \item
	 {\bf  Une lettre perdue de Fermat à Pascal.}
	
	On a la certitude que Fermat a envoyé à Pascal la solution du
	problème des partis avant le 29 juillet 1654, sans doute par
	l'intermédiaire de M. de Carcavi. Nous nous réfèrerons donc
	à une certaine `lettre perdue'.

	 \item
          {\bf Pascal à Fermat, 29 juillet 1654.}
	
	Méthode originale
	pour trouver le parti entre deux joueurs. Calculs
	combinatoires assez poussés.

	\item
          {\bf Pascal à Fermat, 24 août 1654.}

	Pascal critique à tort l'utilisation des combinaisons par Fermat.

	\item
	{\bf Fermat à Pascal, 29 août 1654.}
	
	Politesses. Arithmétique. Remarquons que Fermat
	n'a pas encore reçu la lettre du 24 août et qu'il
	énonce un faux 
	théorème de théorie des nombres resté célèbre%
	\footnote{Je ne
	résiste pas au plaisir 
	idiot  de
	prendre en défaut un génie mathématique : ``Les
	puissances quarrées de 
	2, augmentées de l'unité, sont toujours des nombres
	premiers. [\ldots] C'est une propriété de la vérité de
	laquelle je vous réponds.''}.

	{\bf \item
          Fermat à Pascal, 25 septembre 1654.}

	Réponse à la lettre du 24 août. Fermat justifie brillamment
	sa méthode.

	{\bf \item
          Pascal à Fermat, 27 octobre 1654.}

	Pascal reconnaît la justesse de la méthode de Fermat.

	{\bf \item
  	  Fermat à Carcavi, juin 1656.}

	Fermat donne sans justification les réponses à diverses
	questions, probablement posées par Huygens.

	{\bf \item
 	  Huygens à Carcavi, 6 juillet 1656.}

	Huygens constate que Fermat a trouvé les réponses, et en donne
 	  des démonstrations. Il sollicite l'avis de Pascal et de Fermat
 	  à propos de son traité. 

	{\bf \item
	  Carcavi à Huygens, 28 septembre 1656.}

	Carcavi commente certaines méthodes de Pascal et Fermat.

   \end{itemize}

\item Traités, communications diverses, ouvrages :
   \begin{itemize}
	\item
           \emph{Traité du triangle arithmétique}%
	\footnote{La datation
	du \emph{Traité} pose problème. Comme on l'a vu, Pascal le
	mentionne déjà dans 
	\emph{L'Adresse à l'Académie Parisienne}. Fermat en évoque la
	``onzième conséquence'' dans sa lettre du 29~août~1654. Mais
	rien n'indique que la partie qui nous importe le plus, \emph{Usage
	du triangle arithmétique pour determiner les partis qu'on doit
	faire entre plusieurs joueurs qui jouent en plusieurs
	parties}, fût déjà rédigée à cette époque. Surtout, on
	ignore la forme exacte du \emph{Traité} en 1654 puisque
	selon les \emph{\OE uvres complètes} éditées au Seuil, il  n'a
	été publié qu'en 1665 chez Guillaume 
	Desprez. La même remarque vaut pour \emph{Divers usages du
	triangle arithmétique} et \emph{Combinationes}.},
	de Pascal.

	\item
           \emph{Divers usages du triangle arithmétique}, de Pascal.

	\item
	   \emph{Combinationes}, de Pascal.

	\item
           \emph{Adresse à l'Académie Parisienne}, 1654, de Pascal.

	\item
           \emph{Du calcul dans les jeux de hasard}, 1656\,--\,1657, de
	Huygens. 

     \end{itemize}

\end{itemize}

Dans la suite, nous ne parlerons plus des \emph{Divers usages du triangle
arithmétique}, que nous considèrerons comme une partie du \emph{Traité
du triangle arithmétique}.

\vspace{.4cm}
\subsection*{Ce qui sera dit}

Dans la première partie de ce mémoire, je chercherai à donner une
idée aussi précise que possible des différentes techniques en vigueur
en 1654 pour résoudre le
problème des partis. Nous verrons qu'une lecture de la correspondance entre
Pascal et Fermat fait 
apparaître dès l'origine au moins deux méthodes concurrentes : la méthode
`pas à pas' et la méthode par les combinaisons. Il sera surtout
question de la méthode par les combinaisons, car elle seule pose de
véritables 
problèmes techniques.

Nous nous intéresserons dans la seconde partie aux
fondements de la géométrie du hasard au XVII\ieme\ siècle.
Après m'être expliqué sur ce
que j'entends par fondement, après une étude du \emph{Traité du triangle
arithmétique}, je tenterai de montrer que les intentions, les limites,
l'origine peut-être
de la méthode `pas à pas' s'inscrivent dans la longue tradition
de ces travaux qui, à l'instar des \emph{éléments d'Euclide}, prétendent
fonder les mathématiques sur des bases
claires et solides.

\section{Le savoir-faire de Pascal et de Fermat}

\subsection{La méthode par les combinaisons et la méthode `pas à pas'}
\vspace{.2cm}

\subsubsection*{Quelques lignes étonnantes dans une lettre fameuse}

La lettre de Pascal à Fermat du 29 juillet 1654 est passée à la
postérité. En dehors peut-être de la célèbre marge où Fermat
prétendait avoir 
trouvé la preuve de son grand théorème%
\footnote{Je crois qu'il est inutile de chercher à se procurer une
reproduction de cette marge. On lit en effet dans le Bloc Notes de
Didier Nordon, in \emph{Pour la Science}, Mars 1999, page 5 : ``Le
livre sur lequel Fermat a griffonné serait perdu --- et sa marge avec.''}, 
peu de fragments écrits de la
main d'un 
mathématicien connurent une telle célébrité. D'ou vient pareil renom ?
Tout d'abord cette lettre marque la naissance d'une
branche radicalement nouvelle des mathématiques. Et surtout, ses
premiers paragraphes sont 
peut-être un exemple unique d'un texte mathématique fondateur immédiatement
accessible à quiconque possède des notions élémentaires
d'arithmétique.  Le style en est 
si clair, les arguments si limpides que Pascal emporte la conviction
avec une apparente économie de moyens techniques et
conceptuels. Plutôt que de résumer ce \emph{morceau} (vaine 
entreprise), j'ai préféré 
le reproduire intégralement :

\begin{quotation}

Voici à peu près comme je fais pour savoir la valeur de chacune des
parties, quand deux joueurs jouent, par exemple, en trois parties, et
chacun a mis 32 pistoles au jeu :

Posons que le premier en ait deux et l'autre une; ils jouent
maintenant une partie, dont le sort est tel que, si le premier la
gagne, il gagne tout l'argent qui est au jeu, savoir, 64 pistoles; si
l'autre la gagne, il sont deux parties à deux parties, et par
conséquent, s'ils veulent se séparer, il faut qu'ils retirent chacun
leur mise, savoir, chacun 32 pistoles.

Considérez donc, Monsieur, que si le premier gagne, il lui appartient
64; s'il perd, il lui appartient 32. Donc s'ils veulent ne point
hasarder cette partie et se séparer sans la jouer, le premier doit
dire : ``Je suis sûr d'avoir 32 pistoles, car la perte même me les
donne; mais pour les 32 autres, peut-être je les aurai, peut-être
vous les aurez; le hasard est égal; partageons donc ces 32 pistoles
par la moitié et me donnez, outre cela, mes 32 qui me sont sûres.''
Il aura donc 48 pistoles et l'autre 16.

Posons maintenant que le premier ait deux parties et l'autre point, et
ils commencent à jouer une partie. Le sort de cette partie est tel
que, si le premier la gagne, il tire tout l'argent, 64 pistoles; si
l'autre la gagne, les voilà revenus au cas précédent, auquel le
premier aura deux parties et l'autre une.

Or, nous avons déjà montré qu'en ce cas, il appartient à celui qui a
les deux parties, 48 pistoles : donc, s'ils veulent ne point jouer
cette partie, il doit dire ainsi : ``Si je la gagne, je gagnerai tout,
qui est 64; si je la perds, il m'appartiendra légitimement 48 : donc
donnez-moi les 48 qui me sont certaines au cas même que je perde, et
partageons les 16 autres par la moitié, puisqu'il y a autant de hasard
que vous les gagniez comme moi.'' Ainsi il y aura 48 et 8, qui sont 56
pistoles.

Posons enfin que le premier n'ait qu'\emph{une} partie et l'autre
point. Vous voyez, Monsieur, que, s'ils commencent une partie
nouvelle, le sort en est tel que, si le premier la gagne, il aura deux
parties à point, et partant, par le cas précédent, il appartient 56;
s'il la perd, ils sont partie à partie : donc il lui appartient 32
pistoles. Donc il doit dire : ``Si vous voulez ne la pas jouer,
donnez-moi 32 pistoles qui me sont sûres, et partageons le reste de
56 par la moitié. De 56 ôtez 32, reste 24; partagez donc 24 par la
moitié, prenez-en 12, et moi 12, qui avec 32, font 44.''

\end{quotation}

Voilà le problème des partis brillamment résolu en quelques lignes.
Mais le danger de ce fragment, c'est qu'ébloui par son
éclat, on risque d'oublier qu'il ne s'agit justement que d'un
fragment. Les techniques de Pascal et de Fermat pour calculer les partis sont, dès
1654, assez variées et cette méthode que nous appelons
conventionnellement `pas à pas' n'est qu'un échantillon dont il nous
appartiendra, au-delà de son prestige, de montrer la singularité, tant
dans les travaux de Pascal lui-même, que dans l'histoire générale de la
théorie des probabilités.

Car la méthode `pas à pas', sur l'origine de laquelle nous reviendrons
dans la seconde partie de ce travail, n'est pas la première
imaginée par Pascal pour calculer les
partis. Pascal félicite en effet Fermat 
d'y être parvenu par une autre méthode\footnote{``Votre méthode est très
sûre et est celle qui m'est la première venue'', lettre du 29
juillet 1654.} : ``la méthode qui
procède par les combinaisons''%
\footnote{L'expression est tirée de la lettre du 24 août 1654.}. 
Malheureusement, la solution de Fermat
est aujourd'hui 
perdue. Mais la suite de la lettre de Pascal et les lettres suivantes
nous donnent une  idée
de cette méthode par les combinaisons.

\subsubsection*{Calculs combinatoires}

Dans sa lettre à Fermat du 24 août 1654, Pascal émet des doutes sur
les méthodes de son correspondant. Pour lancer la discussion sur de
bonnes bases, il
commence par résumer le contenu de la lettre perdue, ce qui nous en
donne un aperçu. Voici donc la méthode imaginée indépendamment par
Pascal et Fermat pour résoudre le problème des partis. Je préfère
encore une fois en donner le texte original dans son intégralité~:

\begin{quotation}

Voici comment vous procédez quand il y a \emph{deux} joueurs :

Si deux joueurs, jouant en plusieurs parties, se trouvent en cet état
qu'il manque \emph{deux} parties au premier et \emph{trois} au second,
pour trouver le parti, il faut, dites-vous, voir en combien de parties
le jeu sera décidé absolument.

Il est aisé de supputer que ce sera en \emph{quatre} parties, d'où
vous concluez qu'il faut voir combien quatre parties se combinent
entre deux joueurs et voir combien il y a de combinaisons pour faire
gagner le premier et combien pour le second et partager l'argent
suivant cette proportion. J'eusse eu peine à entendre ce discours-là,
si je ne l'eusse su de moi-même auparavant; aussi vous l'aviez écrit
dans cette pensée. Donc, pour voir combien quatre parties se combinent
entre deux joueurs, il faut imaginer qu'ils jouent avec un dé à deux
faces (puisqu'ils ne sont que deux joueurs), comme à croix et pile, et
qu'ils jettent quatre de ces dés (parce qu'ils jouent en quatre
parties); et maintenant, il faut voir combien ces dés peuvent avoir
d'assiettes différentes. Cela est aisé à supputer : il en peuvent
avoir \emph{seize} qui est le second degré de \emph{quatre},
c'est-à-dire, le quarré. Car figurons-nous qu'une des faces est
marquée $a$, favorable au premier joueur, et l'autre $b$, favorable au
second; donc ces quatre dés peuvent s'asseoir sur une de ces seize
assiettes : 

\vskip.5cm
\[
\begin{array}{|c|c|c|c||c|c|c|c||c|c|c|c||c|c|c|c|} \hline
	a&a&a&a&a&a&a&a&b&b&b&b&b&b&b&b \\
	a&a&a&a&b&b&b&b&a&a&a&a&b&b&b&b \\
	a&a&b&b&a&a&b&b&a&a&b&b&a&a&b&b	\\
	a&b&a&b&a&b&a&b&a&b&a&b&a&b&a&b \\ \hline 
	1&1&1&1&1&1&1&2&1&1&1&2&1&2&2&2 \\ \hline  
\end{array}
\]
\vskip.5cm

\noindent et, parce qu'il manque deux parties au premier joueur, toutes
les faces qui ont deux $a$ le font gagner : donc il en a 11 pour lui;
et parce qu'il y manque trois partie au second, toutes les faces où il
y a trois $b$ le peuvent faire gagner : donc il y en a 5. Donc il
faut qu'ils partagent la somme comme 11 à 5.

\end{quotation}

Un lecteur attentif ne manquera pas de remarquer que le jeu ne se
terminera pas nécessairement en quatre parties. Il se peut très bien que le
joueur 1 gagne ses deux premiers coups et qu'il soit inutile de
poursuivre. Considérer cependant pour la commodité du calcul un jeu se
terminant uniformément en quatre coups est un procédé légitime. 
Cet artifice de 
démonstration, que Pascal appelle le recours à une
``condition  feinte''\footnote{Lettre du 24 août 1654.} fait l'objet
d'une intéressante polémique entre Pascal, Fermat et 
Roberval sur laquelle nous reviendrons par la suite et qui devrait
convaincre les sceptiques. Concentrons-nous 
plutôt, pour le moment, sur les aspects calculatoires de la méthode.
Dans la lettre du 24 août, ils sont tout à fait élémentaires, alors
que dans sa lettre précédente, celle du 29 juillet 1654, Pascal 
fournit des résultats beaucoup plus techniques. Il montre comment,
dans une partie en $n$ coups gagnant, on peut déterminer le parti d'un
joueur menant par une victoire à zéro, ou comme il dit : ``étant donné
tel nombre de parties qu'on voudra, trouver la valeur de la première''.
Le joueur menant par une partie doit, explique Pascal, considérer
comme acquise sa propre mise, mais aussi prendre
une certaine proportion $p'$ de la mise de son adversaire, qui
n'est autre que le quotient du produit des {\small$n-1$} premiers nombres
impairs par celui des {\small{$n-1$}} premiers nombres pairs. Nous
dirions aujourd'hui~:  
\[
    p'=\frac
	{\displaystyle \prod_{i=1}^{n-1} (2i-1)}
	{\displaystyle \prod_{i=1}^{n-1} (2i)}
\]

Cette formule surprend pour qui a un peu l'habitude des
exercices modernes de calcul des probabilités. Comme pouvait l'écrire le pohête Louis
Bouilhet~(il ne parlait pas de mathématiques), cette formule, 
\begin{verse} \ldots\ Elle a cela pour elle,\\
                      Que les sots d'aucun temps n'en ont su faire cas,\\ 
                \protect[\ldots]\ qu'elle est limpide et belle\footnote{Citation
empruntée à \emph{Pauvre Bouilhet} de Henri Raczymow, Gallimard,
1998, page 65, dans le but principal d'essayer l'environnement
{\ttfamily verse}  de \LaTeX.},\\
\end{verse} 

\noindent et une pareille élégance est assez inhabituelle dans une matière où
l'on s'attend plutôt à des entassements d'éléments dénués de
signification apparente.

Pascal parvient à cette formule au prix de calculs assez longs,
et difficiles à suivre, mais dont l'argument central
semble être la valeur $p$ de 
la probabilité%
\footnote{Attention, Pascal n'utilise jamais ce mot, ni aucun
synonyme. Toutefois la lettre du 24 août 1654 nous montre bien qu'il
calculait des ``rapports'' de fait égaux à des probabilités. Nous reviendrons
au \S 1.3 sur le statut de la notion de probabilité au XVII\ieme\ siècle.}
qu'a un joueur de gagner lorsqu'il mène par un à zéro : 

\[	p=\frac{\comb{2(n-1)}{n-1}
		+\comb{2(n-1)}{n}
		+\comb{2(n-1)}{n+1}
		+\cdots
		+\comb{2(n-1)}{2(n-1)}
		}{2^{2(n-1)}}
\]

\vspace{.2cm}

Cette valeur s'obtient par généralisation de la méthode de la lettre du 24
août 1654 : on remarque
que le jeu se termine au plus en $2(n-1)$ manches. Il y a donc
au dénominateur ${2^{2(n-1)}}$ parties différentes possibles, parmi
lesquelles il reste à dénombrer 
celles qui sont gagnées par le joueur 1 (c'est-à-dire celles où il
gagne $(n-1)$ parties, celles où il en gagne $n$ \ldots\ et enfin celle
où il en gagne $2(n-1)$) ce qui donne bien la somme de coefficients du
binôme du numérateur.

Pascal ne fournit pas d'explications complètes sur l'étonnante valeur de $p'$
ni \emph{a~fortiori} sur son lien avec celle de $p$.
Par ailleurs le rapport $p$  n'apparaît même pas explicitement dans sa
lettre. Le
lecteur féru de mathématiques voulant se convaincre que Pascal a
bien calculé $p$ un jour ou l'autre pourra se 
reporter à l'annexe C pour plus de détails.

Comme on le voit, Pascal réalise donc des calculs d'une belle richesse
technique. Loin d'être une simple possibilité abandonnée par la suite, la
méthode des combinaisons est pour lui un objet de recherches
attentives et finalement fructueuses. Mais il a certainement surestimé
l'importance de sa formule : elle n'est, pour autant
que je sache, liée à aucun résultat important, elle ne se généralise
pas, et sa remarquable
symétrie paraît au fond contingente. Pourtant, Carcavi la signale à
Huygens dans sa lettre du  28 septembre 1656 : ``Et le dit S\up{r}
Pascal n'a trouvé la règle 
que lorsqu'un des joueurs a une partie à point ou deux parties à point
(lorsque l'on joue en plusieurs parties), mais il n'a pas la règle
générale''. Cela suggère bien que Pascal s'est fourvoyé, momentanément
au moins, dans la recherche de ``la règle'', d'une belle formule.

Cela dit, comment
est-il parvenu à concevoir ce véritable embryon de calcul des
probabilités fondé d'emblée sur la combinatoire ? Avant d'avancer
quelque hypothèse, il est nécessaire de s'attarder un peu sur le contexte
local des relations entre Pascal et Fermat, mais aussi
tout d'abord sur l'histoire du problème des partis au XVI\ieme\ siècle.

\subsubsection*{Contexte historique}

Dans une certaine mesure, la méthode des combinaisons peut s'inscrire dans la
continuité des tentatives archaïques de résolution du problème des
partis. De Pacioli à
Cardan, l'analyse du problème des partis se raffine.

Pacioli se contente d'un calcul rudimentaire. D'une part il compte le
nombre maximal $n$ de coups
nécessaires à la victoire d'un des joueurs. D'autre part, il compte le
nombre de points de chacun, et partage l'enjeu au prorata de ces points,
sans  faire, curieusement, aucun usage de $n$.
Si la méthode de Pacioli paraît aujourd'hui farfelue, elle n'en est pas
moins 
porteuse d'au moins trois idées importantes :

\begin{itemize}

\item
La première est tout
simplement de se poser le problème des partis. C'est déjà le signe
d'une certaine maturité : après tout, Pacioli aurait pu
dire que celui qui a 50 points a autant de chances de l'emporter
que celui qui en a 20, le hasard se devant de compenser ses
méchants coups par de bonnes fortunes. En fait, si Pacioli s'était
fourvoyé dans de telles croyances, il n'aurait sans doute rien écrit
sur le problème
des partis.

\item
Ensuite vient l'idée que pour résoudre le problème des partis, il faut
faire tôt ou tard un partage au prorata d'un certain nombre de
points. Après tout, on ne fait rien d'autre que cela en calculant
une probabilité, tout l'art étant, comme le fait remarquer Fermat, de
compter des ``hasards égaux''%
\footnote{Lettre du 29 août
1654. Pascal parle aussi, en de nombreux passages de ses textes, de
chances égales, mais dans un contexte très
différent. Il ne s'agit pas pour lui de dénombrer des chances égales
pour calculer ce qu'on appellerait maintenant une probabilité, mais de
justifier le partage en deux de l'enjeu d'une partie.}.

\item
Enfin, et c'est à la fois la plus étonnante et la
moins importante, l'idée de faire intervenir le nombre maximum de parties
indispensable à l'achèvement du jeu, et de se fonder sur ce nombre. La
plus étonnante car comme on va le voir c'est un point crucial de la
méthode par les combinaisons. La moins importante, parce que Pacioli
n'a sûrement pas compris l'intérêt de ce calcul : on a vu
qu'il n'en fait aucun usage.
C'est tout de même un point à méditer, et sur lequel je n'ai pas
grand-chose à dire : 
l'idée de compter le nombre maximal de parties restant à jouer est à la
fois l'idée la plus intéressante et la plus facilement prise en compte par les
prédécesseurs de Pascal et de Fermat. Peut-être parce que ce nombre
est facile à calculer et ne prête pas à polémique. Dans une matière
aussi ardue que le problème des partis, c'est une perche trop belle
pour être négligée. Pacioli s'en saisit et n'en fait rien.

\end{itemize}

Avec Cardan,
l'étape décisive est franchie. Il est le premier qui cherche à
résoudre le problème dans des cas concrets très simples (2 victoires à
1 \ldots), et surtout
il remarque que les victoires engrangées par les joueurs ne comptent pas,
que seul importe le nombre de parties encore à jouer\footnote{E. Coumet
fait en substance les mêmes remarques dans son article. Mais son
optique, plus originale, est de montrer par une fine analyse des
écrits de Cardan que celui-ci est aussi le précurseur de la méthode
`pas à pas'. Je ne reviendrai pas sur ce point, et renvoie donc le
lecteur intéressé à l'article de E. Coumet. \emph{op. cit.}}.

Répétons le, le problème des partis a une histoire avant Pascal. Et
si, comme le dit E. Coumet, Cardan est bien un prédécesseur du
Pascal de la 
méthode `pas à pas', on peut en outre affirmer que les géomètres italiens
de la Renaissance sont des précurseurs du Pascal de la méthode des
combinaisons.
La question de savoir si Pascal ou Fermat ont eu connaissance
de tous ces travaux est
pour autant que je sache ouverte. 

Plus proches de Pascal et de Fermat, les travaux sur les
combinaisons, initiés par Mersenne%
\footnote{Dans \emph{l'Harmonie universelle contenant la théorie et la
pratique de la musique}, Cramoisy, 1636, Mersenne se livre à des
dénombrements. 
Par exemple : combien peut-on écrire de morceaux différents avec huit
notes ? 
}, 
et dont je n'ai malheureusement
pas eu le temps de prendre connaissance, ont déjà en 1654 presque 20
ans. Nous voilà revenus au XVII\ieme\ siècle, une époque où
la vie savante européenne commençait à s'organiser.

\subsubsection*{Le contexte des relations entre Pascal et Fermat}

Les relations entre Pascal et Fermat prennent place dans le cadre des
échanges entre savants. à l'époque, il n'y avait pas de grandes
institutions scientifiques (l'Académie royale des sciences ne sera fondée
par Colbert qu'en 1666) et des hommes tels que le père Mersenne ou M. de
Carcavi aimaient à mettre en rapport les mathématiciens et les
physiciens de toute l'Europe. L'habitude de soumettre des questions à
la sagacité de ses confrères était alors un mode naturel de communication. La
correspondance entre Pascal et Fermat sur le problème des partis
débute précisément par l'entremise
de Carcavi : c'est lui qui pose à Fermat le problème des partis.
Mais de là à croire que l'histoire de ce problème au XVII\ieme\ siècle commence par la
lettre perdue et la lettre du 29 juillet 1654, il y a tout de même
un pas.

Tout d'abord, la lettre perdue n'est pas le premier contact entre Pascal et
Fermat. Il existe une autre lettre%
\footnote{Hélas sa date m'a échappé. Je l'ai
vu dans les \emph{\OE uvres complètes} de Fermat, et je n'ai vraiment
pas le temps d'y retourner. Si je souviens bien, Blaise Pascal devait
alors avoir à peu près 17 ans quand cette lettre a été écrite.}, 
qui n'a rien à voir avec les partis, beaucoup  
plus ancienne, peut-être destinée au père de
Pascal. Cette lettre indique en tout cas que Pascal et Fermat
avaient au moins eu vent l'un de l'autre avant 1654.

 Pour
les questions qui nous occupent, la lettre non datée de Fermat à
Pascal pose un problème plus
délicat et mérite d'être analysée : il y est question d'un problème de jeu
de dés. Un joueur parie
qu'en lançant le dé 8 fois, il fera au moins une fois~6. Après
l'engagement de l'enjeu, le joueur renonce finalement à jouer son
premier coup, et demande à en être dédommagé. Après quoi il renonce
également à jouer son deuxième coup, son troisième, et ainsi de
suite. à chaque 
étape, le joueur reprend donc pour lui un sixième de l'enjeu encore présent, puisqu'il
renonce à avoir une chance sur six de le gagner. Selon Fermat, Pascal
aurait ainsi
calculé la valeur de la $k$\up{ième} partie, qui n'est pas
constante en fonction de $k$ puisqu'à chaque abandon, l'enjeu diminue
de $1/6$. Et Pascal en conclurait qu'un
joueur qui, ayant 
perdu trois parties, renonce à sa quatrième, doit
récupérer cette valeur de la quatrième partie de huit, soit
${125}/{1296}$ du total de la 
mise de départ. Ce en
quoi Fermat le contredit avec juste 
raison : renoncer à la quatrième partie après avoir perdu les trois
premières, c'est exactement renoncer à un sixième de l'enjeu de
départ, car, dit Fermat, ``la somme entière restant dans le jeu, il ne
suit pas 
seulement du principe, mais il est même du sens naturel que chaque
coup doit donner un égal avantage''.

Avant de revenir à nos problèmes de chronologie,
notons au passage que
Pascal semble ici être la victime de l'emploi abusif d'une notion : la valeur
de la k\up{ième} partie de $n$, notion présentant l'inconvénient,
certes inévitable en calculs des probabilités, de difficilement se
transporter d'un problème 
à un autre. Je vois là un travers de philosophe, qui, en quête
d'une certaine universalité, recherche  à toute force et
en toute situation à user d'un
principe général. Sans dire que Pascal à ``l'esprit faux''%
\footnote{Pensée \no 1.},
puisque ``les esprits faux ne sont
jamais ni fins ni géomètres'',
je ne peux m'empêcher ici de lui retourner les
gentillesses qu'il adresse aux ``esprits fins qui ne sont pas
géomètres'', qui ``ne peuvent du tout se
tourner vers les principes de géométrie''%
\footnote{Pensée \no 1.}. 
Il me paraît difficile de savoir si la lettre non datée
est antérieure au reste de la correspondance. Le 29 juillet
1654 il est fait état d'un 
probléme ``des dés'', ce qui peut la replacer dans le contexte de la
correspondance, sans
que l'on puisse alors lui attribuer une date plausible. Fermat y
évoque déjà 
des problèmes subtils, ce qui milite pour une lettre tardive, mais
à l'aide de techniques rudimentaires, ce qui est  au contraire le
signe d'une certaine précocité \ldots\ Heureusement, la suite de la
correspondance nous offre une matière suffisante sans qu'il soit
indispensable de recourir à d'hypothétiques datations.

La lettre du 29 juillet 1654 a visiblement été écrite rapidement :
Pascal y 
mentionne son ``impatience'', pour avouer dans sa lettre suivante
qu'il n'a pu ``ouvrir sa pensée entière touchant les partis''.
Les lacunes dans les démonstrations trouvent peut-être là
leur origine.  L'impression générale se dégageant de la lettre est
celle d'un premier 
contact : les félicitations de Pascal (``J'admire \ldots'') laissent
transparaître un certain étonnement qui nous suggère qu'il
n'avait pas l'habitude de la confrontation avec le génie de Fermat.
Encore une fois, rien ne prouve qu'on ait affaire à la première
lettre concernant la géométrie du hasard.

Sans aller jusqu'à préjuger de l'existence d'une correspondance
perdue entre
Fermat et Pascal, qui expliquerait un peu trop facilement l'irruption
si soudaine et 
simultanée de la vérité dans la question des jeux de hasard, il est
certain qu'en 1654, Pascal n'est pas un franc-tireur s'intéressant
seul à la combinatoire et à la géométrie du
hasard. Il y avait  de nombreuses discussions à ce sujet
entre savants, mais aussi entre esprits curieux qui n'ont pas
forcément laissé leur nom dans l'Histoire : on a vu que le père Mersenne
s'occupait de combinatoire (sans lien avec le hasard), Roberval
discutait avec 
Pascal, Pascal dit à Fermat qu'il communique sa méthode ``à nos
Messieurs''%
\footnote{Lettre du 24 août 1654.}, un certain M. de
Gagnière suggérait à Pascal le moyen de calculer directement les
combinaisons%
\footnote{Mentionné à la fin du
traité \emph{Combinationes}.}, le chevalier de Méré a proposé à Pascal
le problème des partis\footnote{Lettre du 29 juillet 1654. Nous
reviendrons sur ce point dans la conclusion.} \ldots 
\ Si ce bouillonnement a laissé aussi peu de traces, c'est peut-être
qu'avant la percée pascalienne, personne ne pensait que l'étude du
hasard 
mènerait un jour à quelque chose de sérieux. C'est sans doute l'un des  
principaux mérites de Pascal que d'y avoir cru le premier.
La question, que je n'ai
pas résolue, est de savoir si Fermat, du fond de sa province
toulousaine, 
était ou non mêlé à ces savantes mondanités. Sans doute pas : Carcavi,
sceptique ou enthousiaste, a simplement voulu soumettre au grand
géomètre les idées bizarres du
``S\up{r} Pascal''. L'on ne doit pas trop s'étonner que Fermat ait
su résoudre seul le problème des partis : Fermat c'est Fermat
après tout. Tout au plus conviendrons-nous que grâce peut-être à
Mersenne, ce point pouvant sans doute être élucidé, Fermat
savait déjà pas mal de choses sur les combinaisons.

Mais nous n'avons pour l'instant que des informations contextuelles sur 
la place de la méthode par les combinaisons, tant chez Pascal que chez
Fermat. Elle est pourtant l'objet de
jugements surprenants et de controverses au sein même de leur correspondance.

\subsection{Les limites de la méthode par les combinaisons}

Aujourd'hui, la méthode des combinaisons est le moyen naturel, aussi
bien pour l'enseignant, l'ingénieur que le joueur de poker%
\footnote{ En
 ce qui concerne les enseignants et les ingénieurs, n'importe quel ouvrage
 où il est question de probabilités discrètes montre l'importance
 des combinaisons. Pour les joueurs de poker, on consultera avec
 profit \emph{L'illusion du hasard, traité de poker en partie libre},
 par Alexis Beuve, aux éditions Rouge Vif, août 1997.},
de résoudre des problèmes élémentaires de probabilités discrètes. Aussi
peut-il paraître absurde de s'interroger sur les limites de la
méthode des combinaisons. La suite montrera que non.

\subsubsection*{L'aversion de Pascal pour la méthode par les
 combinaisons} 

 Chez Pascal, la méthode par les combinaisons
 est toujours reléguée au second plan
 avec des conséquences parfois énormes. Ainsi, dans le
 \emph{Traité du triangle arithmétique}, Pascal distingue-t-il trois
 méthodes pour ``déterminer les partis qu'on doit faire entre deux joueurs''
 : la méthode `pas à pas', la méthode ``par le moyen du triangle
 arithmétique'', et une fort mystérieuse méthode ``par les
 combinaisons''. Mystérieuse car, étrange négligence, Pascal
 annonce cette méthode\footnote{\emph{Divers usages du Triangle
 arithmétique}, troisième partie, après le ``septième cas''.} sans 
 donner la moindre indication sur son contenu par la suite.

On peut tenter d'expliquer cette lacune.
La méthode ``par les combinaisons'' est
rigoureusement identique à celle ``par le moyen du triangle
arithmétique'', que Pascal donne dans son \emph{Traité}. Les
coefficients inscrits dans ce triangle sont en effet 
ce qu'on appelle aujourd'hui les coefficients du binôme,
des dénombrement de combinaisons. Pascal aurait donc renoncé à écrire cette
platitude que la méthode par les combinaisons n'est autre que celle
par le triangle, au mot près. Et par une
négligence ou une interpolation, l'annonce de la partie
abandonnée est restée là comme un fossile. Cette hypothèse ne me
convainc que partiellement. J'essaierai de la compléter dans la
deuxième partie de ce mémoire. En attendant, Pascal lui-même donne
quelques indications 
sur son rejet de la méthode des combinaisons.

\subsubsection*{``La peine des combinaisons''}

 Selon
Pascal, ``la peine des combinaisons est excessive''%
\footnote{Lettre du 29
juillet 1654, paragraphe 3.}. S'il relègue les combinaisons au
 second plan, ce serait simplement par paresse, par fatigue si l'on
 tient à laisser la religion en dehors de tout ça.
 Voilà une réponse qui nous confronte à un
 nouveau mystère : nous avons vu Pascal manier les
 combinaisons avec virtuosité.
 Je crois qu'il faut voir là un effet de ce qu'Alexandre Koyré
 appelle son ``refus des formules''%
\footnote{\emph{Pascal savant \emph{in} études d'histoire de la
 pensée scientifique} chez Gallimard, collection Tel, 1973.}.  

Le lecteur d'aujourd'hui ne peut qu'être frappé par la façon dont
Pascal présente ses calculs. Par exemple, pour expliquer à Fermat ce
qu'est le produit des huit premiers nombres pairs, après lui avoir
rappelé à toute fin utile que
ces nombres étaient ``2, 4, 6, 8, 12, 14, 16'', il prend soin de
décrire le calcul : 
\begin{quote}

Multipliez les nombres pairs en cette sorte : le premier par le
second, le produit par le troisième, le produit par le quatrième, le
produit par le cinquième etc.

\end{quote}

Je rappelle que Pascal s'adresse au meilleur mathématicien de son
temps, et que la lettre du 29 juillet 1654 d'où est tiré cet
intéressant passage a sans doute été écrite dans une relative
impatience, sans souci de clarté excessive en tout cas. Alors
pourquoi Pascal se lance-t-il dans des explications aussi triviales ? 

à
mon avis, et cela transparaît ailleurs que dans cet exemple grossier,
Pascal n'accorde que peu de signification à une expression telle que `le
produit 
des premiers entiers pairs'. Pour lui une formule n'est que la
retranscription abrégée d'un algorithme. N'oublions pas que Pascal est
aussi l'inventeur de la machine à calculer. Et
``le nombre de façons dont 2 se combine dans 4'', ce que nous notons
$\comb{4}{2}$, n'est pas le résumé concis d'un ensemble de propriétés
algébriques que nous aimons tant, mais plutôt le raccourci trompeur d'un algorithme
complexe. Notons bien que dès le XVII\ieme\ siècle, avec Descartes
(mort en 1650) et l'utilisation d'équations de courbes en
géométrie, le statut des calculs était en train de changer. Plus
encore, en 1656, le traité de Huygens, \emph{Du calcul dans 
les jeux de hasard}, brièvement étudié à l'annexe B de ce mémoire,
montre qu'au  XVII\ieme\ siècle, il n'y avait pas d'obstacle
majeur à une présentation algébrique et relativement moderne de la
géométrie du hasard.
Pascal peut donc paraître un peu en retard sur son temps quant à
sa méfiance envers l'algèbre.

Ainsi, on comprend aisément pourquoi les combinaisons sont 
une peine pour Pascal : si chaque fois qu'il dénombre une
combinaison, il considère l'algorithme plutôt que la formule, si
donc il
garde en tête l'ensemble des calculs effectifs 
que cela sous-tend, comme la lettre du 29 juillet semble l'indiquer,
le calcul le plus élémentaire devient alors étrangement 
ardu. Pire, en lui même il est dénué de sens, tare irrémédiable pour
un philosophe.
Ce dernier point est à mon avis essentiel car Pascal n'est quand
même pas
systématiquement allergique aux formules. Il
semble bien aimer l'expression `la valeur de la $k$\up{ième} partie de
$n$', que l'on rencontre partout et dont il donne des tables dans la
lettre du 29 juillet 1654. Si les usages
l'avaient imposé, il aurait sans doute inventé une notation pour
désigner ce nombre, car on a l'impression --- le calcul de $p'$ en
témoigne --- qu'il en cherche obstinément
les propriétés algébriques.
Sans doute parce que ce nombre a un \emph{sens} lié 
directement au problème des partis.
Et si les combinaisons ont elles aussi un sens, ce n'est pas dans leur
calcul qu'il 
faut le rechercher mais dans ce qu'elles dénombrent, réalité dont le
lien avec le  problème des partis est problèmatique puisqu'il fait
l'objet 
d'une controverse entre trois des plus grands esprits de l'époque.

\subsubsection*{Une controverse entre grands esprits}

Nous avons vu plus haut que Pascal avait des
doutes quant à la 
justesse des raisonnements de Fermat. 
Il n'est pas le seul puisque
Roberval conteste carrément la méthode des combinaisons. Il craint que
l'usage de ``conditions feintes'' ne conduise à  un
``paralogisme''\footnote{Propos rapportés 
par Pascal dans la lettre du 24 août 1654.}.

Plus précisément, ainsi que le rapporte Pascal dans sa lettre du 24
août 1654, il considère ``que c'est à tort que l'on prend l'art de faire
le parti sur la supposition qu'on joue en \emph{quatre} parties, vu
que quand il manque \emph{deux} parties à l'un et \emph{trois} à
l'autre, il n'est pas de nécessité que l'on joue \emph{quatre},
pouvant arriver qu'on n'en jouera que \emph{deux} ou \emph{trois}, ou
à la vérité peut-être \emph{quatre}''.

Pascal rétorque que ``si deux joueurs se trouvant en cet état de
l'hypothèse qu'il manque \emph{deux} parties à l'un et \emph{trois} à
l'autre, conviennent maintenant de gré à gré qu'on joue quatre parties
complètes, c'est-à-dire qu'on jette les quatre dés à deux faces tous à
la fois, [\ldots\ alors] le parti doit être [\ldots] suivant la
multitude des assiettes favorable à chacun''. Roberval admet que ``cela en effet est
démonstratif'', mais ne voit pas en quoi cela justifie le recours à la
feinte.

Pascal explique donc : ``si le premier gagne les deux premières
parties de quatre, et qu'ainsi il ait gagné, refusera-t-il de jouer
encore deux parties, vu que s'il les gagne, il n'a pas mieux gagné, et
s'il les perd, il n'a pas moins gagné ?''. Finalement, ``puisque ces
deux conditions [la naturelle et la feinte] sont égales et
indifférentes, le parti doit être tout pareil en l'une et en
l'autre''. On retiendra que Pascal fait ici preuve d'une belle audace : il
admet implicitement que si deux jeux sont
équivalents pour le sens commun, les calculs qu'on peut faire pour
l'un ou pour l'autre donneront les mêmes résultats. Un jugement
plutôt inné et humain sur l'équivalence de deux situations conduit
donc à une égalité numérique, fantastique découverte à mettre au
crédit de Pascal. Nous verrons toutefois qu'il appartiendra à Fermat
d'en analyser la 
véritable nature mathématique.

Pour le parti à faire entre deux joueurs, Pascal reconnaît donc la
validité de la
méthode par les combinaisons. Avec trois joueurs, il n'en va plus de
même. Dans le cas où ``il manque \emph{une} partie au premier'',
``\emph{deux} au second'' et ``\emph{deux} au troisième'', il pense
qu'on ne peut user de la condition feinte. Car explique-t-il, dans cette
condition, ``deux [joueurs] peuvent atteindre le nombre de leur
parties''. En effet, le jeu se décidera en trois parties au plus. Si
l'on suppose qu'elles sont effectivement jouées, il se peut que le
premier joueur en gagne une, et le deuxième joueur deux, ce qui conduit
selon Pascal à désigner deux gagnants. Dans sa lettre du 25 septembre
1654, Fermat rétorque : ``il semble que vous ne vous souveniez plus
que tout ce qui se fait après que l'un des joueurs a gagné, ne sert
plus de rien''. Quand le premier joueur gagne une partie, ``qu'importe
que le troisième en gagne deux ensuite, puisque, quand il en gagneroit
trente, tout cela seroit superflu ?''

Pour fixer les idées, le tableau suivant, présent dans la lettre du 24 août
1654, permet de calculer le parti. Le dé a trois
faces $a$, $b$, et $c$ respectivement favorables aux joueurs 1, 2 et 3.
Lancé trois fois, il peut donner 27 assiettes. 
Chaque colonne indique une assiette possible, puis 
le vainqueur correspondant. Dans la colonne
\framebox{\begin{tabular}{@{}c@{}}$a$\\$b$\\$b$\\\end{tabular}} , Pascal inscrivait à tort,
en plus du joueur 1, le joueur 2. Fermat, lui, calcule le juste parti
en dénombrant les résultats inscrits dans la dernière ligne :
$17/27$ des enjeux pour le premier joueur, $5/27$ pour les deux autres.

{\[
\begin{array}{|@{}c@{}c@{}c@{}|@{}c@{}c@{}c@{}|@{}c@{}c@{}c@{}|@{}c@{}c@{}c@{}|@{}c@{}c@{}c@{}|@{}c@{}c@{}c@{}|@{}c@{}c@{}c@{}|@{}c@{}c@{}c@{}|@{}c@{}c@{}c@{}|} \hline
	\w{a}&\w{a}&\w{a}&\w{a}&\w{a}&\w{a}&\w{a}&\w{a}&\w{a}&\w{b}&\w{b}&\w{b}&\w{b}&\w{b}&\w{b}&\w{b}&\w{b}&\w{b}&\w{c}&\w{c}&\w{c}&\w{c}&\w{c}&\w{c}&\w{c}&\w{c}&\w{c}  \\
	\w{a}&\w{a}&\w{a}&\w{b}&\w{b}&\w{b}&\w{c}&\w{c}&\w{c}&\w{a}&\w{a}&\w{a}&\w{b}&\w{b}&\w{b}&\w{c}&\w{c}&\w{c}&\w{a}&\w{a}&\w{a}&\w{b}&\w{b}&\w{b}&\w{c}&\w{c}&\w{c}
	\\
	\w{a}&\w{b}&\w{c}&\w{a}&\w{b}&\w{c}&\w{a}&\w{b}&\w{c}&\w{a}&\w{b}&\w{c}&\w{a}&\w{b}&\w{c}&\w{a}&\w{b}&\w{c}&\w{a}&\w{b}&\w{c}&\w{a}&\w{b}&\w{c}&\w{a}&\w{b}&\w{c}
	\\ \hline
	\w{1}&\w{1}&\w{1}&\w{1}&\w{1}&\w{1}&\w{1}&\w{1}&\w{1}&\w{1}&\w{1}&\w{1}&\w{2}&\w{2}&\w{2}&\w{1}&\w{2}&\w{3}&\w{1}&\w{1}&\w{1}&\w{1}&\w{2}&\w{3}&\w{3}&\w{3}&\w{3}
	\\ \hline
\end{array}
\]}

Pascal commet donc une légère erreur dans la justification de l'usage
de la condition feinte. Pour lui, cette condition est légitime si elle ne
conduit qu'à des scores finaux où seul l'un des joueurs détient le
nombre requis de victoires. Fermat fait une analyse
plus fine, dont il s'explique dans la suite de la lettre du 25
septembre 1654.

Car après tout, l'astuce de la condition feinte est-elle réellement
nécessaire au calcul des partis ? Pas le moins du monde. Fermat imagine
une autre solution%
\footnote{Je signale que cette solution est la première que j'ai envisagée
quand je cherchais à résoudre pour moi-même le problème des partis
afin de
m'assurer de la véracité de certains
calculs de Pascal. Cette méthode me semble d'ailleurs être
la première à laquelle pense naturellement 
quelqu'un qui a étudié les probabilités dans un cursus
scolaire actuel. Elle est particulièrement
inefficace pour traiter le cas général d'une partie en $n$ coups
gagnant.}. 
Il commence par préciser que  ``le premier [joueur] peut
gagner, ou 
en une seule partie, ou en 
deux, ou en trois.''. Après quoi il calcule les probabilités de
victoire dans chacun de ces trois cas. Regardons par
exemple comment il calcule la probabilité de victoire en deux parties
du premier joueur :

\begin{quotation} 

Si on en joue deux, il peut gagner de deux façons, ou lorsque le
second joueur gagne la première et lui la seconde, ou lorsque le
troisième gagne la première et lui la seconde. Or, deux dés produisent
9 hasards : ce joueur a donc pour lui $\frac{2}{9}$ des hasards,
lorsqu'on joue deux parties.

\end{quotation}

Fermat calcule donc la probabilité finale en additionnant les trois
probabilités trouvées :

\begin{quotation} 

La somme des hasards qui font gagner ce
premier joueur est par conséquent $\frac{1}{3}$, $\frac{2}{9}$ et
$\frac{2}{27}$, ce qui fait en tout $\frac{17}{27}$.

\end{quotation}

Après cet anodin calcul, qui n'est rien moins que le premier
véritable raisonnement probabiliste de l'histoire des mathématiques, Fermat
clôt la controverse, montrant que depuis le début, il est le seul à
bien comprendre les calculs à l'\oe uvre : il explique que
la condition feinte,
qu'il préfère appeler plus prosaïquement ``l'extension à un
certain nombre de parties'', ``n'est autre que la réduction de
diverses fractions à une même dénomination''. 

Cette controverse fait donc apparaître une hiérarchie : Fermat,
Pascal, puis, bon dernier de la classe, Roberval. On retiendra
simplement la pénétration des vues de Fermat. De tous, il est celui
qui est le plus loin de Paris et des débats qui agitent le monde
savant. Il est pourtant le seul à voir vraiment clair dans les
calculs combinatoires. Mais Pascal qui, de
son propre aveu, ``n'a d'autre avantage sur [Fermat] que celui d'y
avoir beaucoup plus médité''\footnote{Lettre du 24 août 1654.}, était
peut-être handicapé par sa prétendue avance. Qui sait la teneur de
ses discussions avec Méré ou Roberval ? Le problème des partis n'a
peut-être pas été posé d'emblée dans toute sa netteté, des ébauches
de solutions ont pu être proposées avant que l'énoncé ne prenne sa
forme définitive. Ces conversations ont pu laisser de mauvais
alluvions dans l'esprit de Pascal. Ne dit-il pas que parfois ``on se
gâte l'esprit et le sentiment par les conversations''%
\footnote{Pensée \no 6.} ?

\subsection{Comment les fondateurs de la géométrie du hasard n'ont
 pas découvert la notion de probabilité}

Dès 1654, les grands mathématiciens français se préoccupaient donc de
problèmes de probabilités avec un vocabulaire et des techniques très
proches de ce qui se pratique aujourd'hui encore. Mais paradoxalement,
le grand absent 
dans leur débat est justement la notion de probabilité d'un
événement. Elle n'apparaît que dans la méthode sans condition feinte
de Fermat (``\ldots $\frac{2}{9}$ des hasards''). 

De fait, Pascal et Fermat calculent des probabilités, mais sans voir,
sans dire en tout cas, que le fameux rapport des cas favorables et des
cas totaux est, en un sens assez précis%
\footnote{Ce sens se dégage des définitions et des théorèmes de la
théorie moderne des probabilités, de la loi faible des grands nombres
notamment.},
une mesure des chances de gagner des joueurs. 
Et si Pascal évoque ``la fraction qui exprime la valeur de la première
partie de huit''\footnote{Lettre du 29 juillet 1654.}, c'est pour
ajouter immédiatement après : ``c'est-à-dire que si on joue chacun le
nombre de pistoles exprimé par [le dénominateur], il en appartiendra
sur l'argent de l'autre le nombre exprimé par [le numérateur]''. 

Il semble donc que
la probabilité n'est à l'époque qu'un coefficient nécessaire au partage des enjeux.
On conviendra que la méthode `pas à pas' est directement
fondée sur la notion d'espérance. Mais la méthode des combinaisons,
d'où l'on s'attendrait à voir émerger la notion de probabilité, 
semble elle aussi paradoxalement imprégnée de
l'idée d'espérance. Ainsi, s'il fallait donner un nom moderne à la
géométrie du hasard selon 
Pascal, ce serait bien mieux `théorie de l'espérance' que `théorie des
probabilités'. 

En guise d'explication, j'indiquerai simplement que
la notion de 
probabilité n'a \emph{a priori}
que peu d'intérêt une fois sortie de son contexte. Sans support théorique,
associer le nombre $1/2$ 
au jeu de pile ou face n'est qu'une paraphrase abstraite et numérique de l'idée
naturelle de
chances égales. Là n'est pas la richesse de la probabilité qui
provient plutôt de
son lien avec le concept d'événement : la probabilité de la
disjonction ou de la conjonction
de deux événements se calcule fort bien moyennant les définitons
\emph{ad hoc}. Si la notion de probabilité a
fini par s'imposer au fondement de la théorie, c'est parce qu'elle est
conservée par les opérations d'une
algèbre des événements. Fermat s'engage dans cette voie, comme on l'a
vu, en calculant sans feinte le parti dans le cas de trois joueurs,
mais sans donner suite. Pourquoi ? 

En premier lieu parce que dans le cas d'un jeu à $n$ joueurs, sa
méthode probabiliste ne permet pas facilement
d'exprimer le parti par une formule ramassée. Pour autant que mes
tentatives n'aient pas été trop maladroites, la mise au même
dénominateur qui se fait sans souffrance par le moyen de la condition
feinte est inextricable quand on se trouve nez à nez avec des
monstruosités du type :

\[
\frac{1}{2^2}\comb{2}{2}+\frac{1}{2^3}\comb{3}{2}+\cdots+\frac{1}{2^{k+2}}\comb{k+1}{2}
\]%
Cette formule est directement issue de l'application de la
méthode de Fermat. Au premier joueur, il manque 2 parties, à l'autre
k. Le jeu peut donc se terminer en 2, 3, \dots, ou $k+2$ coups. Dans
chacun des cas, on dénombre les cas favorables au premier joueur, soit
par exemple dans le cas où la partie se terminerait en $k+1$ coups,
$\comb{k+1}{2}$. On divise ce nombre par le nombre total de cas
envisageables dans une partie en $k+2$ coups : $2^{k+2}$.

Mais la complexité des calculs n'est peut-être pas ce qui a détourné
Fermat des méthodes probabilistes. En fait, la géométrie
du hasard dans son ensemble ne l'intéresse sans doute pas. En 1654, sa
correspondance montre qu'il 
s'occupe de
théorie des nombres%
\footnote{La théorie des nombres n'est évidemment pas totalement déconnectée
de la combinatoire. Les combinaisons peuvent notamment être utilisées
pour prouver ce qu'on appelle le Petit Théorème de Fermat : Si
$a\geq 1$
est un entier et $p$
un entier premier, alors $a^{p-1} \equiv 1$ modulo $p$. J'ignore
toutefois si Fermat a utilisé la preuve par les combinaisons
pour démontrer sa célèbre congruence. Cela expliquerait en tout cas son
aisance avec les combinaisons, sans pour autant contredire son absence
d'enthousiasme pour la géométrie du hasard.}, 
et les questions de Pascal devaient n'être pour lui
que des problèmes plaisants, certes subtils, mais indignes de la
vraie mathématique.

Bref, Fermat, en pur géomètre, ne s'intéresse pas à un
problème qui en 1654 ne semblait pas promis aux développements
\emph{mathématiques} 
éclatants qu'il a eu par la suite. à propos, Fermat n'est-il pas le
premier 
(un des premiers disons, car je
n'en sais rien !)
grand mathématicien à n'avoir été que mathématicien, inaugurant ainsi
une division des savoirs contre laquelle il est de bon ton aujourd'hui
de vociférer ? Si bien qu'on se demande si Pascal n'était pas ironique, ou
plus probablement 
maladroit, quand dans le flot croissant de compliments qu'il adresse à
Fermat, il finit par lui écrire qu'il fait ``peu de différence entre
un homme qui n'est que géomètre et un habile
artisan''\footnote{Lettre du 10 août 1660, tirée  des \emph{\OE
uvres complètes} de Pascal, au Seuil.}, 
pour immédiatement, bien sûr, l'exonérer de ce travers.
Pascal, en bon
philosophe, est au contraire stupéfait de voir ``la fortune incertaine [\ldots] si
bien maîtrisée par l'équité du
calcul''\footnote{\emph{Combinationes}.}. On comprend donc qu'il ait
souhaité aller au-delà des questions purement techniques et
calculatoires de la géométrie du hasard, sur des chemins que nous
allons maintenant explorer en son agréable compagnie.

\newpage
\section{Les fondements de la géométrie du hasard au XVII\ieme\ siècle}

\subsection{Position du problème}

Parler du fondement de la théorie des probabilités au XVII\ieme\ siècle
peut paraître anachronique. On considère généralement que la théorie 
des probabilités est établie sur des bases théoriques
solides
depuis les travaux de Kolmogorov dans les années 1930. Mais si l'on oublie un 
peu le sens actuel, formaliste pour dire vite, du fondement d'une théorie
mathématique, on doit bien admettre que dès le XVII\ieme\ siècle, la
nouvelle ``géométrie du 
hasard'', suspecte dans son appellation même, méritait bien, dans un
sens ou un autre, d'être fondée
sur quelque principe au moins, à défaut d'axiomes au sens où nous
l'entendons depuis Hilbert. Mon point de vue n'est cependant pas de
parler des fondements de la géométrie du hasard dans un sens faible
(ordonner en un tout plus ou
moins cohérent les techniques éparses du moment), mais bien dans le sens
fort, celui qui convient aux travaux sur les fondements de l'analyse
du XIX\ieme\ siècle par exemple.

Plutôt que d'attaquer de front le problème des fondements de la
géométrie du hasard, je me propose de revenir sur l'énigme de la
méthode `pas à
pas'. Au premier abord, il semble en effet qu'elle n'ait pas de
nécessité mathématique, ni 
même historique. On imagine fort bien l'histoire
générale des
probabilités sans la méthode `pas à pas'. Je donne peut-être ici
l'impression de sombrer dans  une erreur fort classique. Les
mathématiques ont en 
effet l'étrange particularité qu'on les imagine facilement (et à tort)
sans histoire du 
tout, et \emph{a fortiori} sans Pascal. Les \emph{éléments de
mathématique} de Bourbaki ne prennent-ils pas les mathématiques ``à
leur début''%
\footnote{Je ne dis pas que Bourbaki considère que les mathématiques
n'ont pas d'histoire. Mais les traités de Bourbaki commencent par un ``Mode
d'emploi'' dont la première phrase est invariablement : ``Ce traité
prend les mathématiques à leur début, et donne des démonstrations
complètes. Sa lecture ne suppose donc en principe aucune connaissance
mathématique particulière''. Il faut voir là une preuve de l'humour de
Bourbaki, humour attesté par Raymond Queneau dans son essai
\emph{Bourbaki et les mathématiques de demain, \emph{in} Bords} chez
Hermann, 1978.},
dans un sens qui n'a rien d'historique ?
Imagine-t-on un traité de botanique qui prendrait la botanique à son
début ?

Spéculer sur la nécessité historique de tel ou tel développement des
mathématiques  est
donc souvent, je l'admets, tout à fait douteux. Toutefois, je
maintiens 
que la méthode `pas à pas' n'a pas de nécessité historique
visible : Fermat semble se 
soucier comme d'une guigne des méthodes récurrentes de Pascal. Il
n'y fait pas reférence une seule fois; pour lui les combinaisons ne
sont pas une peine, ne posent pas de problème majeur. Il est
vrai qu'on l'a vu les manier avec une sûreté et une profondeur de vue
inégalées à son époque.  
Bref, entre les mains de Fermat, on conçoit une géométrie du hasard
plus conforme dès sa jeunesse à ce qu'elle est
aujourd'hui. J'insiste : je ne me contente pas d'imaginer ce qu'aurait
pu être tel ou tel développement des mathématiques qui n'a pas eu lieu. La
correspondance montre que Fermat maîtrisait la technique d'une théorie 
élémentaire des probabilités correspondant grosso modo au
programme de nos chères terminales scientifiques. Hélas ou
heureusement, qu'importe, Fermat n'a pas écrit de traité de théorie
des jeux. 

La méthode `pas à pas' est donc singulière à tous égards : elle n'a pas
été une étape nécessaire du développement de la théorie des
probabilités, puisque
la méthode des combinaisons, ancêtre apparent de nos conceptions
actuelles, lui est antérieure. Pis, elle n'a été l'étape de rien du tout,
n'étant pas le point de départ d'une éventuelle
théorie concurrente ou complémentaire de la théorie classique des
probabilités (du moins, cette théorie n'a-t-elle  pas encore vu le jour). 
Afin de percer à jour les intuitions \emph{a posteriori} mystérieuses
de Pascal, nous
allons donc nous 
lancer dans l'étude de son \emph{Traité du triangle arithmétique}, qui
présente l'avantage, par rapport à la correspondance, de donner non pas
forcément le point de vue intime de Pascal, mais ce qu'il
en voulait laisser voir à ses contemporains, à la postérité, à
qui enfin se préoccupe de la Vérité.

\subsection{Le \emph{Traité du triangle arithmétique}}
\vspace{.2cm}
\subsubsection*{Vue d'ensemble}

Sur la question du hasard et des mathématiques, le traité
sérieux le plus ancien dont nous disposons est celui de Pascal, 
avec toutes les réserves de rigueur quant à la datation de sa forme
définitive. La méthode `pas à pas' y est developpée dans toute sa
généralité. Voici pour fixer les idées le plan du \emph{Traité} :

\begin{itemize}

\item
Pascal commence par définir ce qu'est le triangle
arithmétique, aujourd'hui communément appelé
triangle de Pascal. Il donne des propriétés du système de coordonnées
servant à en repérer les cellules. La figure 1 ci-dessous fournit une
reproduction aussi fidèle que possible du triangle arithmétique,
également reproduit sur l'encart détachable.

\item
Il donne ensuite des ``conséquences'' concernant les égalités 
rencontrées dans les cellules du triangle. Les démonstrations de ces
égalités ne mobilisent que des outils très rudimentaires.

\item
Il consacre la suite de son \emph{Traité} aux
proportions qui se rencontrent dans ces cellules. Là, les démonstrations
sont beaucoup plus techniques, et le raisonnement par récurrence,
chose inhabituelle pour l'époque, joue un rôle majeur.

\item
Enfin, il pose et résout
le problème de calculer le contenu d'une cellule du triangle de manière
directe. 

\item
Il fournit un traité annexe, \emph{Divers usages du triangle
arithmétique}, qui contient quatre parties :

\begin{itemize}

\item
``Usage du triangle arithmétique pour les ordres numériques'', qui
donne une application du triangle à la recherche des nombres
polygonaux.

\item
``Usage du triangle arithmétique pour les combinaisons'', où Pascal
montre que les nombres du triangle arithmétique sont des dénombrements
de combinaisons.

\item
``Usage du triangle arithmétique pour déterminer les partis qu'on doit
faire entre deux joueurs qui jouent en plusieurs parties'', qui
présente en toute généralité la méthode `pas à pas'.

\item
``Usage du triangle arithmétique pour trouver les puissances des
binômes et des apotômes'', qui présente l'importante formule du
binôme. 

\end{itemize}

\end{itemize}

\begin{figure}
\begin{center}

\begin{tabular}{r|c|c|c|c|c|c|c|c|c|c|}  
&{1}&{2}&{3}&{4}&{5}&{6}&{7}&{8}&{9}&{10}\\ \hline
1&\z{\y{\mathtt{G}}{1}}&\z{\y{\sigma}{1}}&\z{\y{\pi}{1}}
&\z{\y{\lambda}{1}}  &\z{\y{\mu}{1}}  &\z{\y{\delta}{1}}
&\z{\y{\zeta}{1}} &\z{1} &\z{1}&\z{1} \\ \cline{1-11}
2&\z{\y{\varphi}{1}}&\z{\y{\psi}{2}}&\z{\y{\theta}{3}}
&\z{\y{\mathtt{R}}{4}}  &\z{\y{\mathtt{S}}{5}}  &\z{\y{\mathtt{N}}{6}}
&\z{7} &\z{8} &\z{{9}} \\ \cline{1-10}
3&\z{\y{\mathtt{A}}{1}}&\z{\y{\mathtt{B}}{3}}&\z{\y{\mathtt{C}}{6}}
&\z{\y{\omega}{10}}  &\z{\y{\xi}{15}}  &\z{{21}}  &\z{{28}} &\z{{36}}
\\ \cline{1-9}
4&\z{\y{\mathtt{D}}{1}}&\z{\y{\mathtt{E}}{4}}&\z{\y{\mathtt{F}}{10}}
&\z{\y{\rho}{20}}  &\z{\y{\gamma}{35}}  &\z{{56}}  &\z{{84}} \\ \cline{1-8}
5&\z{\y{\mathtt{H}}{1}}&\z{\y{\mathtt{M}}{5}}&\z{\y{\mathtt{K}}{15}}
&\z{{35}}  &\z{{70}}  &\z{{126}}  \\ \cline{1-7}
6&\z{\y{\mathtt{P}}{1}}&\z{\y{\mathtt{Q}}{6}}&\z{{21}}  &\z{{56}}
&\z{{126}}  \\ \cline{1-6}
7&\z{\y{\mathtt{V}}{1}}&\z{{7}}&\z{{28}}  &\z{{84}}  \\ \cline{1-5}
8&\z{{1}}&\z{{8}}&\z{{36}}  \\ \cline{1-4}
9&\z{{1}}&\z{{9}}\\ \cline{1-3}
10&\z{{1}}\\ \cline{1-2}
\end{tabular}
\caption{Le triangle arithmétique.}
\end{center}
\end{figure}

La première ligne et la première colonne ne font pas, à
proprement parler, partie du triangle. Elles ne servent qu'à en
repérer les cellules. Pascal
appelle ``exposants'' les nombres qu'elles
contiennent.
Rappelons comment sont définies les cellules du triangle arithmétique :
on place le nombre 1 à la cellule G, (le ``générateur''). Chaque cellule
du triangle est alors définie comme étant la
somme de la cellule qui se trouve à sa gauche, et de celle qui se trouve
au-dessus. Par exemple, $C=\theta+B$. Les colonnes sont
appelées ``rangs perpendiculaires'' et les lignes ``rangs
parallèles''.  Chaque cellule a donc deux exposants (que nous
appellerions plus volontiers coordonnées) : son ``exposant parallèle'', et
son ``exposant perpendiculaire''. Les lignes obliques joignant des
exposants égaux (celle qui 
contient 1, 3, 3 et 1 par exemple) sont appelées les ``bases'' du
triangle. La ligne diagonale (qui contient 1, 2, 6, 20, 70 \dots) est
appelée la ``dividente''.

Une lecture superficielle du \emph{Traité} pose déjà des problèmes amusants. Il
faudra qu'un jour un amateur d'études kabalistiques se penche
sur la question du nom donné à chaque cellule du triangle, qui repose sur
une nomenclature aujourd'hui tombée en désuétude%
\footnote{ On
m'objectera que cette pratique subsiste encore dans le choix les sigles en
vigueur à l'administration
du M.E.N.. Je rétorque que le recours aux
lettres grecques n'y est pas encore trop répandu.} : 
chaque cellule est
désignée par une lettre latine ou grecque, le rapport entre la
position de la cellule et sa lettre échappant à toute logique
raisonnable. 

Les amateurs de l'histoire du nombre zéro --- il y en a
--- 
prêteront une attention particulière à la démonstration de la
Conséquence première, où, pour prouver que le premier rang
perpendiculaire du triangle ne contient que des 1, Pascal fait
remarquer que l'application du principe de génération des
cellules revient, pour les cellules de ce rang, à additionner un zéro 
à leur cellule d'au-dessus.

Mais ce qui frappe le plus aujourd'hui, c'est la manière
géométrique dont Pascal définit le triangle. Il mène explicitement des
parallèles, des perpendiculaires~\ldots toutes choses
aujourd'hui 
secondaires, les relations algébriques entre les cellules étant
privilégiées. Gardons-nous de ne voir là qu'un des
archaïsmes pittoresques qui parsèment les textes
d'époque. La suite prétend au contraire montrer que la géométrie du hasard
de Pascal était réellement imprégnée de géométrie. 

Tout au long du \emph{Traité}, Pascal convient implicitement de confondre
abusivement une cellule et 
le nombre qu'elle contient. à la Conséquence neuvième, il est
même prouvé qu'une base
est ``égale'' à une 
somme d'autres bases, alors qu'en toute rigueur, nous dirions
aujourd'hui que l'égalité concerne les nombres contenus dans les
bases. On reconnaît là une habitude de la géométrie ancienne
présente déjà dans les \emph{éléments} d'Euclide, qui consiste à dire
par exemple que 
deux triangles sont égaux si, dirions-nous aujourd'hui, on les peut
superposer modulo une isométrie du plan, ou que deux angles sont égaux
s'ils ont même mesure \ldots\ Cela suggère que
Pascal espérait véritablement construire une \emph{géométrie} du hasard, où
les nombres
du triangle joueraient un rôle comparable aux
notions de surface ou d'angle. D'ailleurs, quand Pascal fait
intervenir le nombre zéro, dans la preuve de la
Conséquence première, il pense
peut-être plus à une quantité métrique (surface, angle nul \ldots)
qu'à l'élément neutre de l'addition, autrement dit dans un jargon
moins furieusement tendance, 
ce qui reste quand on n'additionne rien du tout. Mais c'est là
un point accessoire, il y en a de plus probants :

Au
début du \emph{Traité}, Pascal remarque que la somme des deux exposants
d'une cellule ne dépend que de la base où elle se trouve. Ce résultat évident,
``plutôt compris que démontré'', repose selon Pascal sur la
division ``en un pareil nombre de partie'' des deux côtés orthogonaux
du triangle. Je ne vois pas bien le lien démonstratif entre ces deux
propositions, mais je vois un penchant certain à se référer à des
raisonnements métriques à la Thalès. En tout cas, si l'on devait
réécrire tout cela aujourd'hui, on ne tenterait en aucun cas de se ramener à
des propriétés du triangle, mais plutôt à de laborieuses considérations
sur des incrémentations d'indices.

Sans tirer de conséquences trop audacieuses de ce qui vient
d'être avancé, on peut 
néanmoins conclure que Pascal, grand lecteur d'Euclide, cherchait à
imiter, inconsciemment 
peut-être, ce qui se faisait encore de mieux en matière de rigueur
et de démonstration à son époque : la géométrie euclidienne. Plus
avant, nous en rencontrerons un indice plus profond.

Le style agréable de Pascal mérite aussi une mention, car il risque de
tromper un lecteur pressé. Ses textes mathématiques sont en effet écrits
dans une langue remarquable, avec tant d'aisance qu'on peut avoir
la fausse impression qu'il ne démontre pas grand chose, qu'il se
contente de décrire avec élégance quelque vérité que lui aurait suggérée
son intuition. Cette impression est notablement accentuée par l'usage
original qu'il fait du \emph{nombre générique}. Là où un auteur
moderne désignerait les cellules du triangle par des
variables, Pascal préfère utiliser 
des exemples particuliers, ce qu'on considère généralement comme
monstrueux,  les mathématiques n'ayant que faire de cas particuliers.
Mais il le fait avec une telle habilité que la
rigueur n'en pâtit pas. Regardons par exemple comment Pascal prouve
la Conséquence seconde de son traité :

\vspace{.5cm}

\newpage
\begin{quotation}

{\begin{center} {\noindent \textsc{conséquence seconde}} \end{center}}

\emph{
En tout Triangle arithmétique, chaque cellule est égale à la somme de
toutes celles du rang parallèle précédent, comprise depuis son rang
perpendiculaire jusqu'au premier inclusivement}

Soit une cellule quelconque $\omega$ : je dis qu'elle est égale à
$\mathtt{R}+\theta+\psi+\varphi$, qui sont celles du rang parallèle supérieur
depuis le rang perpendiculaire de $\omega$ jusqu'au premier rang
perpendiculaire.

Cela est évident par la seule interprétation des cellules par celles
d'où elles sont formées.

Car $\omega$ égale $\mathtt{R}+\mathtt{C}$.

\hspace{2.5cm}$\overbrace{\theta+\mathtt{B}}$

\hspace{2.8cm}$\overbrace{\psi+\mathtt{A}}$

\hspace{3.4cm}$\overbrace{\varphi}$

car $\mathtt{A}$ et $\varphi$ sont égaux entre eux par la [conséquence]
précédente. 

Donc $\omega$ égale $\mathtt{R}+\theta+\psi+\varphi$.
\end{quotation}

\vspace{.5cm}

évidemment, à nos yeux, la cellule $\omega$ n'a rien de ``quelconque'',
puisque $\omega$ n'est pas une variable mais une cellule particulière
du triangle.  Pourtant,
on voit bien que le raisonnement fait avec $\omega$ peut être
recommencé avec n'importe quelle cellule du triangle. C'est tout
l'art du nombre générique : écrire un raisonnement sur un cas
particulier, mais en prenant bien garde que le raisonnement puisse
être reproduit mécaniquement sur les autres cas. Dans ce cas précis,
on voit à quelle économie de notation conduit
 un usage habile du nombre générique.

\subsubsection*{Le raisonnement par récurrence}

À partir de la Conséquence douzième, le \emph{Traité} prend une
tournure différente. D'abord, Pascal ne s'intéresse plus aux 
``égalités qui se rencontrent dans le triangle'', mais aux
``proportions''. Il s'agira par exemple de calculer des rapports de cellules
adjacentes du triangle. Cette distinction, aujourd'hui un
peu  spécieuse, se retrouve curieusement dans les \emph{éléments}
d'Euclide, ce qui mérite quelque éclaircissement :  le théorème de
Thalès est le premier qui dans
l'ordre de démonstration des résultats de la géométrie classique
fasse obligatoirement intervenir la notion de nombre réel. Le premier,
car le
théorème de Pythagore par exemple peut s'exprimer en termes de
surfaces et se prouver par des
`couper-coller'%
\footnote{Techniquement, dans les \emph{éléments},
le théorème de Pythagore n'est pas prouvé par des `couper-coller' mais
par la constance de la surface d'un triangle lorsque qu'un de ses
sommets se déplace le long de la droite parallèle au côté
opposé. Mais la démonstration de cette constance fait quant à elle intervenir
des `couper-coller' justifiés par des cas d'égalité des triangles,
qui sont donc en dernière analyse l'argument
utilisé.}. 
En effet, le théorème de Thalès se rapporte
aux proportions, chose difficilement réductible à des notions de pure
géométrie. Une théorie des nombres réels est donc indispensable, et est
effectivement présente au Livre cinquième des \emph{éléments}. 

Par une étrange analogie, Pascal marque une rupture au moment de
s'attaquer aux proportions. Cette rupture ne concerne pas seulement
la nature des résultats énoncés, mais aussi leur démonstration : le
raisonnement par récurrence y tient une place inhabituelle pour
l'époque. Pour nous en faire une idée précise, admirons un peu cette
Conséquence douzième : 

\begin{quotation}

{\begin{center} {\noindent\textsc{Conséquence Douzième}} \end{center}}

\emph{
En tout Triangle arithmétique, deux cellules contiguës étant dans
une même base, la supérieure est à l'inférieure comme la multitude
des cellules depuis la supérieure jusqu'au haut de la base à la
multitude de celles depuis l'inférieure jusqu'en bas inclusivement}

\end{quotation}

Avec ${E}$ et ${C}$, par exemple, qui sont sont bien
contiguës dans une même base, Pascal dit que 
${E}$ est à ${C}$ comme 2 à 3, ce que nous
dirions aujourd'hui $\frac{{E}}{{C}}=\frac{2}{3}$.
2~``parce qu'il y a deux cellules [dans la base] depuis ${E}$
jusqu'en bas, savoir ${E}$, ${H}$'' et~3~``parce qu'il y a
trois cellules depuis ${C}$
jusqu'en haut, savoir ${C}$, ${R}$, $\mu$''. 

La démonstration commence par énoncer le principe de la démonstration
par récurrence. Premier ``lemme'' : la proportion est évidente dans la seconde base car
``$\psi$ est à $\sigma$ comme 1 à 1''. D'autre part, second lemme, si elle ``se
trouve dans une base quelconque, elle se trouvera nécessairement dans
la base suivante.'' Avant de démonter ce dernier lemme, Pascal
explique : 

\begin{quotation}

D'où il se voit qu'elle [la proportion] est nécessairement dans toutes
les bases : car elle est dans la seconde base par le premier lemme;
donc par le second elle est dans la troisième base, donc dans la
quatrième et à l'infini.

\end{quotation}

On voit que le raisonnement par récurrence n'est pas chez Pascal un moyen
technique imaginé pour faire l'économie du concept d'infini actuel
dans des démonstrations ayant ``une infinité de cas'', rôle que d'autres
mathématiciens ont pu en d'autres temps lui faire tenir.
Pascal n'hésite pas à dire ``et
ainsi à l'infini''; après quoi, il se garde bien de mettre des points
de suspension. D'ailleurs, dans la preuve de la Conséquence huitième,
écrite sans raisonnement par récurrence, Pascal termine aussi son
argumentation par un courageux ``et ainsi à l'infini''. Pour ce qui
est de son
utilisation en mathématiques, l'infini ne l'effrayait pas. 

Mais
revenons à la conséquence douzième car la démonstration de son second
lemme est intéressante. étant donné le temps que j'ai passé à
essayé de la comprendre, j'ai l'impression que, faute de
notations modernes, elle manque de clarté. L'usage systématique du
nombre générique 
montre ici ses limites. Quand Pascal parle de 2, de 3 ou de 4, on ne
sait plus du tout d'où viennent ces nombres, et on a du mal à saisir le
caractère démonstratif du discours. Voici donc une version
modernisée, mais fidèle dans ses enchaînements de la preuve de Pascal : 

\vspace{.2cm}

Pour une cellule $X$, notons $V_X$ sa valeur, $L_X$ son rang dans sa
base compté à partir de la gauche et $R_X$ son rang dans
sa base compté à partir de la droite. Par exemple, $L_B=2$ car en
partant de la gauche, dans la base de B on rencontre 2~cellules : D et
B. Et $R_B=3$ car en partant dans la base de B on rencontre 3~cellules
: $\lambda$, $\theta$ et $B$. Nous parviendrons à nos fins en
trois étapes :

\begin{enumerate}

\item 

Par hypothèse de récurrence, $\frac{V_D}{V_B}=\frac{L_D}{R_B}$

\vspace{.1cm}

Donc,  $\frac{V_D}{V_B}+1=\frac{L_D}{R_B}+1$ et
$\frac{V_D+V_B}{V_B}=\frac{L_D+R_B}{R_B}$.

\vspace{.1cm}

Or, $V_D+V_B=V_E$ par construction du triangle arithmétique. Donc,

\[\frac{V_E}{V_B}=\frac{L_D+R_B}{R_B}\]

\item

De même, par hypothèse de récurrence, $\frac{V_B}{V_\theta}=\frac{L_B}{R_\theta}$

\vspace{.1cm}

Donc,  $\frac{V_\theta}{V_B}+1=\frac{R_\theta}{L_B}+1$ et
$\frac{V_\theta+V_B}{V_B}=\frac{R_\theta+L_B}{L_B}$.

\vspace{.1cm}

Or, $V_\theta+V_B=V_C$ par construction du triangle arithmétique. Donc,

\[\frac{V_C}{V_B}=\frac{R_\theta+L_B}{L_B}\]

\item 

Pour conclure, utilisons les propriétés du système de coordonnées
repérant les cellules : $L_D+R_B=R_\theta+L_B$, car ces deux sommes
sont égales au nombre $S$ de cellules de la base contenant les
cellules contiguës $D$, $B$ et
$\theta$.

D'autre part, $R_B=R_C$ car C est à droite de $B$, et $L_B=L_E$ car
$E$ est en dessous de $B$. 

Finalement, en reportant ces égalités dans les égalités trouvées aux
points 1. et 2. on trouve :

\begin{center}
$\frac{V_E}{V_B}=\frac{S}{R_C}$ \hspace{.5cm} et \hspace{.5cm} $\frac{V_C}{V_B}=\frac{S}{L_E}$
\end{center}

Et en divisant ces deux égalités, on obtient bien par la simplification
que Pascal appelle ``la proportion troublée'' : 

\[\frac{V_E}{V_C}=\frac{L_E}{R_C}\]

\end{enumerate}

Si j'inflige ce kroupnik\footnote{\'Etouffe-chrétien.}, au lecteur, ce n'est
pas par manie professionnelle
d'assommer mes semblables avec des calculs, mais pour montrer encore une fois
le caractère géométrique des méthodes de Pascal. Il me semble
naturel de suivre cette démonstration le triangle sous les yeux. à
chaque étape, l'\oe il ou l'index en parcourt les cellules pour
persuader la raison. 
Si $R_B=R_C$, n'est-ce pas parce que $\lambda B C\mu$ est un
parallélogramme ?
N'a-t-on pas l'impression que si $V_B$ se
simplifie à la fin du calcul, c'est au fond parce que $B$ se trouve au
\emph{milieu} de $D$ et $\theta$ ? La ``proportion troublée'' n'est-t-elle
pas le pendant algébrique d'une sorte de relation de Chasles ? Mais
étant donné le caractère
subjectif de tout ceci, je laisse le lecteur en juger par lui-même \ldots

À partir de la Conséquence douzième, Pascal peut égrainer
toute la
série de proportions qui en découle, pour finalement aboutir à
``l'accomplissement de [son] traité'' : le calcul direct de la valeur
d'une cellule connaissant seulement son exposant parallèle et son
exposant perpendiculaire.

J'espère avoir été convaincant : si le nom de l'ancêtre de la théorie
des probabilités fut ``géométrie du hasard'', c'est bien dans le sens
précis du mot géométrie qu'il faut l'entendre.

\subsubsection*{Le problème des partis dans le \emph{Traité}}

Dans la troisième partie des \emph{Divers usages du triangle
arithmétique}, Pascal explique comment faire le parti entre
deux joueurs à l'aide du triangle arithmétique. Sa méthode repose sur
les deux principes suivants :

\begin{itemize}

\item 

``Si un des joueurs se trouve en telle condition que, quoi qu'il
arrive, une certaine somme lui doit appartenir en cas de perte et de gain,
sans que la hasard la lui puisse ôter, il n'en doit faire aucun
parti, mais la prendre entière comme assurée parce que le parti devant
être proportionné au hasard, puisqu'il n'y a nul hasard de perdre,
il doit tout retirer son parti.''

\item

``Si deux joueurs se trouvent en telle condition que, si l'un gagne, il
lui appartiendra une certaine somme, et s'il perd, elle appartiendra à
l'autre; si le jeu est de pur hasard et qu'il y ait autant de hasard
pour l'un que pour l'autre et par conséquent non plus de raison de
gagner pour l'un que pour l'autre, s'ils veulent se séparer sans
jouer, et prendre ce qui leur appartient légitimement, le parti est
qu'ils séparent la somme qui est au hasard par la moitié, et que
chacun prenne la sienne.''

\end{itemize}

En voyant ces principes, on se demande aujourd'hui spontanément
s'ils sont une définition de ce qu'est le parti, ou s'ils sont des lois
en régissant le calcul. Bien sûr, chez Pascal, on se trouve dans le
second cas, les concepts mathématiques n'ayant pas à son époque
le statut arbitraire que nous leur connaissons. Mais en grand
géomètre, Pascal manipule ses principes exactement comme s'il
s'agissait de définitions. C'est-à-dire que presque rien d'extérieur
à eux ne vient étayer les démonstrations. Je dis `presque' car, on le
voit bien, les  principes présentent une lacune : dans n'importe
quel jeu, du
moins dans le genre de jeu auquel Pascal s'intéresse, nul ne peut être
assuré d'un quelconque gain avant la victoire. Comme peut-il alors y
avoir une ``somme assurée'' ? Bien sûr, se poser pareille
question dénote
une compréhension erronée des principes de Pascal, qui ne doivent être
utilisés qu'à bon escient lors de démonstrations
récurrentes, comme on l'a vu dans la méthode `pas à pas'. Ce qui
finalement nous
conduit bien à déceler une lacune :
 c'est l'usage bien compris des principes qui
déterminera seul les endroits précis où l'on est en droit de les
utiliser. 

Si on devait réécrire selon nos critères de rigueur les
principes de Pascal,  on préciserait
donc sans doute que la ``certaine somme'' qui  ``doit appartenir [au
joueur] en cas
de perte et de gain, sans que la hasard la lui puisse ôter'' doit
s'entendre comme une certaine somme provenant d'un calcul qui
utilise le premier principe. On parvient donc à éradiquer la lacune au
prix d'une apparente inconsistance logique : le premier principe fait
usage dans sa propre définition du premier principe. En fait il n'y là
aucune inconsistance : quiconque à
l'habitude de l'informatique, de la logique mathématique ou de la
linguistique formelle dira simplement, chacun selon ses usages, 
qu'on est en présence d'une
définition récursive, inductive%
\footnote{Le mot `induction' est utilisé ici dans le sens qu'on lui donne en
logique mathématique, à ne pas confondre avec l'``opération mentale
qui consiste à remonter des faits à la loi'' chère aux sciences
naturelles. Voir par exemple  \emph{Logique mathématique} par
René Cori et Daniel Lascar, tome 1, chez Masson, 1994, page 20.},
ou encore générative, du parti à faire entre deux joueurs%
\footnote{Gardons-nous
cependant de dire que Pascal est le premier dans l'histoire à faire
usage de ce type de définitions. On peut
lire dans \emph{Compilation, principes, techniques et outils}, 1989,
chez Interéditions, ouvrage connu des informaticiens sous le
sobriquet de `dragon rouge' : ``L'érudit [indien] Panini imagine une
notation syntaxique pour spécifier la grammaire du Sanscrit'', au
III\ieme\ siècle avant J.-C.. Ces notations syntaxiques, appelées
syntaxes BNF (Backus Naur Form) en informatique,
sont des ensembles de règles inductives de spécification d'une syntaxe.
On les retrouve en linguistique sous le nom de grammaires génératives,
inventées par Noam Chomsky.
Les adeptes de l'érudit P\={a}nini peuvent à ce qu'il paraît consulter
\emph{Panini-Backus form suggested}, Ingerman, 1967, Communication ACM,
10-3, p.137.}.
Le parti à un moment $t_1$ du jeu dépend certes du parti associé à un
autre moment $t_2$ du jeu, mais $t_2$ est antérieur à $t_1$. Et par la
remarque ---
étrangement semblable au principe de descente infini de Fermat ---
qu'une partie ne peut durer depuis un temps infini, puisque s'il est
possible d'aller vers l'infini, il est bien sûr impossible d'en
venir, on voit que les applications successives du principe mèneront à
la détermination du juste parti en un nombre fini d'étapes. 

Cette longue digression peut laisser la désagréable impression d'une
réinterprétation \emph{a posteriori} des textes de Pascal. Toutefois,
il est certain 
que les deux principes \emph{sont} inductifs, ce qui est à mettre au
crédit du génie de Pascal. Ce n'est pas tout : ils
frappent par leur  actualité. De la loterie nationale
aux marchés boursiers dérivés en passant par les assurances, notre
époque en est toute imprégnée. Dans notre société, on peut acheter
ou vendre des biens qui ne sont rien d'autre que de l'espérance
mathématique.  Pascal est le premier à avoir découvert (ou inventé ?)
l'existence de ce bien immatériel. 

Après avoir énoncé ses deux principes, Pascal expose la manière dont 
ils s'appliquent à la méthode `pas à pas', qui apparaît alors
pour ce qu'elle est vraiment : une démonstration non pas par
récurrence mais par induction. Si ces passages du
\emph{Traité} n'ont 
pas connu la célébrité des quelques lignes de la lettre du 29 juillet
1654, c'est parce qu'à mon avis ils sont moins clairs. Pascal commet
la maladresse de commencer son induction par le bas%
\footnote{Précaution élémentaire quand on critique quelqu'un comme
Pascal, prendre comme caution un grand 
mathématicien : ``The feeling that this sort of
treatment [les récurrences] adds to the precision of an inductive
argument is much too common and is responsible for the introduction of
many irrelevancies in the literature''. \emph{On numbers and games},
Jonh H. Conway, Academic Press Inc. (London) Ltd, 1976, page 64.
}, 
c'est-à-dire en traitant d'abord les cas les plus simples : 

\begin{quotation}

{\it
\begin{center}\noindent Premier cas \end{center}

Si à un des joueurs il ne manque aucune partie, et à l'autre
quelques-unes, la somme entière appartient au premier.} Car il l'a
gagnée [\ldots]

{\it
\begin{center}\noindent Second cas \end{center}

Si à un des joueurs il manque une partie et à l'autre une }[\ldots]

{\it
\begin{center}\noindent Troisième cas\end{center}

Si à un des joueurs il manque une partie et à l'autre deux }[\ldots]

\end{quotation}

\[
\vdots
\]

Ensuite Pascal peut démontrer par récurrence l'utilité des cellules du
triangle arithmétique pour le problème des partis, utilité qui ne
surprendra personne : ces cellules servent à calculer des
combinaisons%
 \footnote{Pascal montre, à la deuxième partie des
 \emph{Divers usages du triangle artihmétique}, 
 que les cellules du
 triangle dénombrent des combinaisons, grâce à ce qu'on appelle
 aujourd'hui la formule de Pascal :
 $\comb{n}{p}=\comb{n-1}{p-1}+\comb{n-1}{p}$.},
et les combinaisons servent à calculer les partis. Mais
cela, Pascal ne l'a pas dit dans son traité, bien qu'il le sût
\ldots\  
La démonstration du lien entre le triangle et le problème des partis
est donc nécessairement un peu alambiquée, et fait évidemment appel au
raisonnement par récurrence :

Pascal explique d'abord que pour trouver le parti entre deux joueurs
à qui il manque $n$ parties, il faut prendre la $n$\up{ième} base du
triangle, la séparer en deux blocs de cellules contiguës, de sorte
que chaque bloc ait autant de cellules qu'il y a de victoires manquantes
pour chaque joueur. Le partage se fait au prorata de ces deux
nombres. Par exemple, s'il manque deux parties au premier joueur et 
trois au deuxième, le parti est de $\frac{1+4+6}{16}$ et
$\frac{4+1}{16}$, comme la méthode `pas à pas' le montrait d'ailleurs
dans la lettre du 29 juillet 1654. La récurrence se fait sur $n$. Le
passage de base en base du triangle repose sur les deux principes 
et les propriétés des cellules prouvées au début du \emph{Traité}.

\subsection{Nouvel éclairage sur la méthode `pas à pas'}

Au terme de notre exploration des textes de Pascal et de Fermat,
nous pouvons enfin reprendre l'analyse des fameuses lignes de la
lettre du 29 juillet 1654.

\subsubsection*{La place centrale d'une méthode désuète}

Revenons tout d'abord sur les motivations de Pascal. Avec sa géométrie
 du hasard, il est conscient
 d'inventer une théorie nouvelle qui met en
jeu des conceptions profondément ancrées chez l'homme : le hasard, la
providence \ldots Et si on imagine la 
très docte Académie parisienne admettant la supériorité de la méthode
des combinaisons sur toute autre considération (Tartaglia \ldots), si on
l'imagine encore prêtant un statut de vérité géométrique au calcul
des partis, on conçoit  mal en revanche que Pascal ait pu se contenter
du concensus des savants de l'époque pour affermir ses vues. Il lui
fallait quelque chose de solide, de sûr.

Or, du côté des combinaisons, on est dans le règne des polémiques, des
justifications hasardeuses, où Pascal lui-même perd pied. La notion
de probabilité, esquissée par Fermat, devait avoir pour Pascal le
statut d'un de ces abus de langage commodes, qui permettent
d'appréhender un problème et, pourquoi pas, d'en trouver sans faille la
solution, mais en aucun cas de le résoudre en conscience et en vérité.
Car Pascal est un géomètre, et son souci, nous
l'avons montré, est d'intégrer sa géométrie du hasard à la géométrie
tout court : le \emph{Traité du triangle arithmétique} est d'inspiration
euclidienne et, par un véritable tour de force (l'utilisation de
l'induction et de la récurrence), Pascal parvient à asseoir le calcul
des partis sur deux
principes naturels, des axiomes serait-on tenté de dire, desquels
découlent  tous les
calculs.

Bref, en un sens Pascal a bien fondé une géométrie du hasard. 
Nous ne prétendrons pas comparer l'approche de Pascal avec celle de
Kolmogorov 
par exemple, ce qui serait cette fois franchement anachronique. Mais
je pense que le \emph{Traité du triangle arithmétique} garde une certaine
actualité dans la mesure où les principes sur lesquels il fonde
les probabilités sont des principes d'un sens commun dont nos théories
ne se soucient plus guère, ce qui fait d'ailleurs leur force
moderne. Les 
axiomes de Kolmogorov n'ont en effet rien d'axiomes au sens de vérités
évidentes. La théorie moderne des probabilités est abstraite. Elle a
toutefois un point d'ancrage sur le sens commun : le point de vue
fréquentiste. Dire qu'un événement a une probabilité $p$ de se
produire implique qu'une série de $n$ épreuves fera se produire cet
événement avec une fréquence qui tend à se rapprocher de
$p$ quand $n$ tend vers l'infini. Ce résultat se démontre en théorie
des probabilités, mais ce qui nous intéresse ici, le fait qu'il dise
réellement quelque chose sur le monde, n'est plus aujourd'hui du
ressort des mathématiques. Le seul avantage du point de vue
fréquentiste reste
donc qu'il peut se vérifier expérimentalement pour de grandes valeurs
de $n$. J'aimerais bien savoir
ce que Pascal aurait pensé de ce point de vue, lui qui, ne
l'oublions pas, a écrit prophétiquement que l'infini ``est cependant un 
nombre''\footnote{Pensée \no 233.}.

En attendant, on s'explique enfin pourquoi Pascal ne
parle pas de la méthode par les
combinaisons dans
son traité. Sûrement pas parce c'est une méthode qu'il
maîtrise mal : après sa correspondance avec Fermat, il détient un
savoir lui permettant d'éviter ses erreurs passées. C'est
plutôt parce que selon Pascal, cette méthode fonde ses preuves sur un
jugement vague, sans qu'on arrive à les réduire à des
principes évidents. On me rétorquera que Pascal, grâce à
son \emph{Traité} pouvait très bien justifier la méthode par les
combinaisons, les cellules du triangle faisant le lien avec le
problème des partis. Mais de la sorte, il aurait vidé la méthode par
les combinaisons de sa substance, qui est de se fonder sur des
dénombrements. Car par le moyen du triangle arithmétique, il aurait
certainement pu montrer que le nombre qui est égal à une certaine combinaison
sert à calculer les partis. Mais il aurait été bien en peine d'expliquer
dans les canons de la rigueur le lien entre la nature de la
combinaison (un dénombrement) et le
juste parti, il n'aurait pu montrer par quel principe sa méthode
diffère de celle de Pacioli. Aujourd'hui encore, selon moi, ce
principe reste mystérieux, plus compris que démontré. Les
constructions de Pascal, quant à elles, offrent une approche
raisonnable du hasard, mathématiquement dépassée, mais satisfaisante pour
l'esprit. 

La méthode `pas à pas' ou, comme Pascal l'écrivait à Fermat, ``mon
autre méthode universelle, à 
qui rien n'échappe et qui porte sa démonstration avec
soi''\footnote{Lettre du 24 août 1654.}, est donc la brillante synthèse d'un
travail de grande ampleur. Sans en avoir l'air, elle incorpore tout
l'outillage technique et conceptuel du \emph{Traité du triangle
arithmétique}. On peut lui reprocher  d'être
particulière : on voit mal comment l'appliquer à des problèmes
complexes de boules colorées tirées d'urnes avec ou sans remise. Je lui
reprocherais surtout de n'être pas naturelle. La première idée, la
bonne idée, qui
vient pour résoudre le problème des partis, dès les essais de Pacioli et
aujourd'hui encore, est bien de compter les chances de chacun. 
Voyons cela comme un indice supplémentaire : c'est un défaut habituel
des travaux de fondement des mathématiques que de s'embarrasser de
constructions aussi peu
naturelles que l'identification d'un nombre réel à une
coupure de l'ensemble des nombres rationnels.

\subsubsection*{Origine de la méthode `pas à pas'}

Quand Pascal a-t-il imaginé la méthode `pas à pas' ? Nous n'avons pas de
documents précis sur la question et en sommes donc réduits à des
spéculations. On peut toutefois prendre appui sur deux indices :

\begin{itemize}

\item
À propos de la méthode des combinaisons, Pascal écrit à Fermat :
``Votre 
méthode est très sûre et est celle qui m'est la première venue à la
pensée dans cette recherche''\footnote{Lettre du 29 juillet 1654.}

\item
À propos du calcul direct d'une cellule du triangle par les
exposants --- ``l'accomplissement'' de son traité --- Pascal écrit :
``M. de Gagnières me 
communiqua lui-même cette excellente solution et me proposa même
d'en chercher la 
démonstration ;~j'admirai le problème, mais effrayé par la difficulté,
je pensai qu'il convenait d'en laisser la démonstration à son auteur;~%
cependant, grâce au triangle arithmétique, une voie aisée me fut
ouverte pour y parvenir''\footnote{\emph{Combinationes}, paragraphe
antépénultième.}.

\end{itemize}

Le premier indice montre, nous le savions déjà, que Pascal a trouvé la
méthode par les combinaisons avant la méthode `pas à pas'. 
Quant au deuxième, il montre que Pascal savait déjà une partie de ce
qu'il voulait prouver quand il a commencé à s'intéresser au
triangle, que le triangle n'a pas été pour lui un moyen de découverte
mais de démonstration. On peut donc supputer que la méthode `pas à
pas' est à la fois postérieure aux découvertes combinatoires de Pascal
et à la méthode par la combinaison. 

Pascal aurait donc trouvé la méthode `pas à pas' en ayant tous les
outils en main : il savait que
les partis sont liés aux combinaisons, il savait comment manipuler
ces dernières à l'aide du triangle arithmétique. Se dessine alors un
scénario possible : il
cherche à fonder la géométrie du hasard. Il voit dans les méthodes
par récurrence du triangle le moyen technique de parvenir à ses
fins. L'analyse fine du lien 
entre triangle et combinaisons lui suggère alors une adaptation
directe des méthodes
du triangle au problème des partis. Finalement il découvre plutôt
qu'il n'invente 
ses deux principes après coup, comme des justifications des passages
de base en base utilisés lors de la récurrence prouvant le lien entre
le triangle et les partis. Les récurrences et
l'induction sont alors le fruit d'une reflexion sur la nature de ces
justifications. La méthode `pas à
pas' vient en dernier, synthèse de ce travail, et non point de
départ. 

Nous savons déjà que Pascal envisage la méthode `pas à pas' comme un
moyen de fonder sa géométrie du hasard. Cela aurait pu lui apparaître
après qu'il l'a inventée. Si  nos hypothèses sont exactes, ce ne
serait pas ainsi que les choses se passèrent : la méthode `pas à pas'
serait bien une découverte issue de recherches sur les  fondements de la
géométrie du hasard.

\newpage
\section{Conclusion}

Les débuts de la géométrie du hasard nous ont permis de côtoyer deux
des plus grands esprits du XVII\ieme\ siècle. Ceux qui pensent
que les mathématiciens n'ont pas de style devraient s'intéresser à ces
deux là : 

Fermat, le pur mathématicien qui, avec son
intuition presque innée des probabilités, ne se trompe jamais. 
Qui pourtant ne voit pas en cette nouvelle géométrie du hasard 
la grande découverte  qu'elle réprésente pour nous, 
à la manière de ces Vickings qui accostèrent au
Canada sans découvrir l'Amérique. 

Et Pascal, dont les textes
traduisent une véritable obsession de plier la nouvelle géométrie du 
hasard aux canons de l'authentique rigueur mathématique dans ce
qu'elle a de plus noble : éclairer l'esprit. Ce programme
insensé, il n'a pu en venir à bout qu'au prix de
l'invention de nouvelles démarches  : l'induction et l'usage
systématique du raisonnement par
récurrence nous sont légués pour toujours. Pourtant, les idées de
Pascal  sur la manière de
fonder cette science sont absentes de la théorie des probabilités
telle qu'on la pratique. 

C'est que,  comme le disait Ludwig Wittgenstein,
``les problèmes \emph{mathématiques} de ce que l'on appelle les
fondements sont aussi peu pour nous au fondement des mathématiques que
le rocher peint supporte le château peint.''%
\footnote{Citation
empruntée à \emph{Le philosophe et le réel, entretiens avec
Jean-Jacques Rosat}, p.18, par Jacques Bouveresse, Chez Hachette, 1998.}
L'induction et le raisonnement par récurrence, problèmes
\emph{mathématiques} 
du fondement de la géométrie du hasard, ne la fondent en rien, et sont
aujourd'hui remplacés par d'autres techniques, d'autre concepts. Et ce qui
devrait vraiment fonder la géométrie du hasard, les deux principes de
Pascal,  ne sont
plus  du ressort de nos mathématiques.

On regrettera peut-être que je n'ai pas parlé d'éventuelles  conditions
extérieures aux développement des mathématiques, qui auraient
favorisées les travaux de Pascal et Fermat. Sur ce point, je pense
qu'une recherche sur le chevalier de Méré, qui a publié des ouvrages
juridiques, serait intéressante. Après tout, c'est lui qui pose le
problème des partis à Pascal. Sur l'état général des
idées au XVII\ieme\ siècle, on pourra se reporter à
un article de Norbert Meunier, \emph{L'émergence d'une mathématique du
probable au XVII\ieme\ siècle}\footnote{\emph{Revue d'histoire des
mathématiques}, 2 (1996), p. 119--147.}.

Mais on regrettera peut-être surtout qu'il n'ait pas été question
de philosophie dans tout ce qui précède. Curieusement, Pascal
lui-même parle 
 très peu de sa géométrie du hasard dans son \oe uvre
philosophique. En toute franchise je n'ai pas lu tout Pascal. 
Des recherches rapides mais efficaces m'ont pourtant convaincu qu'en
dehors du \emph{Pari} et des 
textes proprement mathématiques, il n'y a rien dans son \oe uvre sur
la géométrie du 
hasard. Peut-être est-ce à cause de l'expérience mystique du \emph{Mémorial}, qui
le 23 novembre 1654 va changer la vie de Pascal et l'engager dans
d'autres aventures. Peut-être aussi parce que
le \emph{Traité du triangle arithmétique} est finalement à ranger au
sein même de
l'\oe uvre philosophique de Pascal.   

\appendix
\newpage

\addcontentsline{toc}{section}{Annexes :}

\renewcommand{\thesection}{\Alph{section}}

\section{Le \emph{Pari} de Pascal}

Je ne m'étendrai pas sur la question du \emph{Pari} de Pascal%
\footnote{Je m'appuie sur le texte de la Pensée \no 233 dans l'édition
du Livre de poche des \emph{Pensées}.}, 
liée à des
questions théologiques et philosophiques sur lesquelles je suis à peu
près incompétent. Le texte du \emph{Pari} est tiré des \emph{Pensées} de Pascal,
qui sont les brouillons d'un ouvrage sur lequel il travaillait :
\emph{Apologie pour la Vérité de la Religion Chrétienne}. On ne sait
pas exactement quand il commença à composer cette \oe uvre, mais en 1658 il
classa ses notes ``en diverses
liasses''\footnote{Information tirée de l'introduction aux
\emph{Pensées} des \emph{\OE uvres complètes} éditées au Seuil.}.
La grande expérience mystique de Pascal, relatée dans son
\emph{Mémorial}\footnote{Page 249 dans l'édition Livre de poche des
\emph{Pensées}.}, date quant à elle, par une coïncidence dont je laisse le
lecteur libre de tirer ce qu'il lui plaira, du 23 novembre 1654, soit moins d'un
mois après la dernière lettre de Pascal à Fermat touchant au problème
des partis, sa lettre du 27 octobre 1654.

Dans l'optique de ce mémoire, le texte du \emph{Pari}
présente la curiosité d'être l'unique endroit où, à la manière de
Fermat,  Pascal fait
intervenir des notions de probabilités relativement conformes à ce
qu'on rencontre aujourd'hui. Pascal considère une probabilité tendant
vers 0, (les chances selon un incroyant d'accéder à la vie éternelle
par la piété) et un 
gain qui est lui infini (cette vie éternelle justement) : 

\begin{quotation}
Mais il y
a une infinité de vie et de bonheur. Et cela
étant, quand il y aurait une infinité de hasards dont un seul serait
pour vous,  vous auriez encore raison de gager un pour avoir deux, et
vous agiriez de mauvais sens, étant obligé à jouer, de refuser de
jouer une vie contre trois à un jeu où d'une infinité de hasard il y
en a un pour vous, s'il y avait une infinité de vie infiniment
heureuse à gagner. Mais il y a ici une infinité de vie infiniment
heureuse à gagner, un hasard de gain contre un nombre fini de hasards de
perte, et ce que vous jouez est fini. Cela ôte tout parti : partout
où est l'infini, et où il n'y a pas infinité de hasard de
perte contre celui de gain, il n'y a point à balancer, il faut tout
donner. 
\end{quotation}

Mathématiquement, cela se tient à peu près, exception faite de
``l'infinité de hasard dont un seul serait pour vous'', qui donne
l'idée d'une probabilité nulle d'accéder au gain infini, et conduit
donc à une
espérance de gain finale délicate à calculer~:~$0\times\infty$.
Toutefois, Pascal ne retombe pas dans cette erreur par la suite,
puisqu'il prend soin de préciser que son argumentation est valable partout ``où il
n'y a pas infinité de hasard de perte'', ce qui exclut la probabilité
nulle --- à mettre sur le compte de l'ardeur d'une plume un peu
leste. Rappelons-nous 
que le texte du \emph{Pari} n'est qu'un brouillon.

Toujours en rapport avec notre sujet, une discussion passionnante sur
l'infini précède le passage du 
\emph{Pari}. C'est une étrange coïncidence vu le rôle que l'infini
était appelé à jouer par la suite en théorie des probabilités
(point de vue fréquentiste, $\sigma$-additivité, mesure de Lebesgue
\ldots). Pascal lui-même a certainement dû méditer sur l'infini
lorsqu'il mettait au point son principe de démonstration par
récurrence et ses méthodes par induction. Mais l'infini dont il est
question dans les \emph{Pensées} est introduit dans un cadre un peu
extérieur aux probabilités : il s'agit surtout de justifier
mathématiquement l'idée d'un gain infini à un jeu, un jeu dans lequel
nous sommes tous embarqués.

\section{Le traité de Huygens}

En 1656, Huygens publie un petit traité, \emph{Du calcul dans 
les jeux de hasard}, qui fera date puisque le
\emph{Traité du triangle arithmétique} de Pascal n'a été publié qu'en
1665%
\footnote{Voir l'article de Norbert Meunier pour une étude plus
complète. \emph{op. cit.}
}.
Cet ouvrage ingénieux part de la même idée que celui de Pascal :
donner un fondement rigoureux à la géométrie du hasard naissante. Mais
il a moins d'ampleur. Destiné à indiquer rapidement comment résoudre
quelques problèmes, il se présente donc différemment. 
Chez Huygens, la proposition centrale est la suivante : 

\vspace{.5cm}
\begin{quotation}

{\begin{center} {\noindent\textsc{Proposition III}} \end{center}}

Avoir $p$ chances d'obtenir $a$ et $q$ chances d'obtenir $b$, les
chances étant équivalentes, me vaut $\frac{pa+qb}{p+q}$.

\end{quotation}
\vspace{.5cm}

On est manifestement plus proche que chez Pascal des idées
actuelles. On comprend bien aujourd'hui l'intérêt de la proposition
de Huygens, que l'on perçoit comme une définition de l'espérance
mathématique. Mais Huygens, dans la lignée de Pascal, a souhaité 
prouver son résultat. Il donne donc un principe : ``dans un jeu, la
chance qu'on a de gagner quelque chose a une valeur telle que si l'on
possède cette valeur on peut se procurer la même chance  par un jeu
équitable''. Ce principe me semble quelque peu confus, et l'usage
qu'en fait Huygens dans sa Proposition I permet de l'éclairer :

\vspace{.5cm}
\begin{quotation}

{\begin{center}{\noindent\textsc{Proposition I}} 

Avoir des chances égales d'obtenir $a$ ou $b$ me vaut $\frac{a+b}{2}$.
\end{center}}

Afin de non seulement démontrer cette règle, mais aussi de la
découvrir, appelons $x$ la valeur de ma chance. Il faut donc que,
possédant $x$, je puisse me procurer de nouveau la même chance par
un jeu équitable. Supposons que ce jeu soit le suivant. Je joue $x$
contre une autre personne, dont l'enjeu est également $x$. Il est
convenu que celui qui gagne donnera $a$ à celui qui perd. Ce jeu est
équitable, et il appert que j'ai ainsi une chance égale d'avoir $a$ en
perdant, ou $2x-a$ en gagnant le jeu; car dans ce dernier cas
j'obtiens l'enjeu $2x$, duquel je dois donner $a$ à l'autre joueur. Si
$2x-a$ était égal à $b$, j'aurais donc une chance égale d'avoir $a$ ou
d'avoir $b$. Je pose donc $2x-a=b$, d'où je tire la valeur de ma
chance $x=\frac{a+b}{2}$. La preuve en est aisée. En effet, possédant
$\frac{a+b}{2}$, je puis hasarder cette somme contre un autre joueur
qui mettra également $\frac{a+b}{2}$, et convenir avec lui que le
gagnant donnera $a$ à l'autre. J'aurai de force une chance égale
d'avoir $a$ si je perds, ou $b$ si je gagne; car dans ce dernier cas
j'obtiens l'enjeu $a+b$ et je lui en donne $a$.

\end{quotation}
\vspace{.5cm}

Cette démonstration, d'un abord malaisé, procède par analyse et
synthèse. Recherchant un jeu satisfaisant le principe, on trouve
d'abord une
équation, dont la solution donne le paramètre pour un bon candidat. On
vérifie ensuite que ce candidat est bien un jeu satisfaisant le
principe. 

On peut se demander à quoi sert la partie `analyse' de la
preuve. Le candidat-jeu est plus facile à trouver par l'intuition que
par le calcul rocambolesque de Huygens, et d'après le principe, son
existence suffit à prouver le théorème. à mon avis, si Huygens s'est
quand même lancé dans une analyse délicate, c'est parce qu'il
sentait que son principe avait une faiblesse logique. Si l'on peut se
contenter de l'existence d'un certain jeu auxiliaire pour trouver une
valeur, qu'est-ce qui prouve alors qu'on ne trouvera pas un autre jeu
auxiliaire, avec une autre valeur ? On voit bien qu'il faut ajouter
quelque part
une condition d'unicité, sinon du jeu auxiliaire, au moins de sa
valeur. Or l'ajout d'une telle condition rend le principe soit
inapplicable, si l'on impose vraiment  une impossible unicité, soit
difficile à 
formuler, si pour sauver l'unicité, on cherche à se restreindre à
une  certaine classe de
jeux. Huygens coupe donc la poire en deux : il ne change rien à son
principe, mais il résout quand même une équation pour faire quelque
chose qui ressemble à une preuve d'unicité.

Nous voilà renseignés sur la manière de fonder le calcul des
probabilités chez Huygens. Cette méthode présente quand même
l'intérêt d'utiliser des jeux auxiliaires dans les preuves. Cette
idée était déjà en germe chez Cardan\footnote{Cardan calculait des partis
(de faux partis évidemment) en cherchant la mise que paierait un autre
joueur pour 
rentrer sans avantage ni dommage dans le jeu. Voir E. Coumet,
\emph{op. cit.}}, et pour la première fois, un jeu devient un objet
mathématique au même titre qu'un nombre ou une figure, et non plus
un simple objet d'étude.

Mais, au fait, comment Huygens procède-t-il  pour trouver le
parti entre deux joueurs à l'aide de la seule Proposition III, qui ne
peut être utilisée qu'une fois les ``chances'' calculées ?
Allez voir par vous même, et vous apprendrez qu'il utilise la
méthode `pas à pas', avec un peu moins d'élégance que Pascal !

On retiendra du mémoire de Huygens son style très
actuel. Les calculs et les preuves sont présentés à quelques détails
près comme dans bien des livres de mathématiques d'aujourd'hui. On
voit là que même selon les critères du XVII\ieme\ siècle, Pascal,
visiblement nostalgique de la géométrie grecque, travaillait un
peu à l'ancienne.

\section{Précisions techniques sur la lettre du 29 juillet 1654}

La lettre du 29 juillet 1654 présente des difficultés de lecture assez
importantes, dues à l'absence totale de démonstration des résultats
obtenus par la méthode des combinaisons. Voici, en langage
mathématique actuel, une démonstration des résultats de
Pascal, respectant l'ordre exacte dans lequel Pascal dit les avoir prouvés, et 
devant donc être assez proche de la démonstration originale (perdue,
ou peut-être jamais mise au net).

Pour qui serait choqué par les formules et les notations modernes, que
Pascal refusait comme le dit Alexandre Koyré, je précise
que ces formules et ces notations peuvent tout à fait 
être exprimées dans le langage de Pascal. Une lecture de la lettre du
29 juillet 1654 permet de s'en convaincre. En outre, essayer d'écrire
ce qui suit `à la manière de Pascal' serait un exercice de
style certes intéressant, mais aussi idiot que courir le
marathon de New-York en sandales.

\subsection*{Remarques préliminaires}

Je m'appuie sur les \emph{\OE uvres complètes} de Pascal publiées au
Seuil en 1988. Je me réfère aux numéros de pages, et à une numérotation
personnelle mais logique des paragraphes.

Comme nous allons le voir, une lecture précise du texte de Pascal est
gênée par l'habitude actuelle de formuler tous les problèmes de
théorie des jeux en termes de probabilité de gagner et d'espérance de gain.

Pour y voir plus clair, commençons par poser exactement le problème auquel s'attaque
Pascal. Deux joueurs engagent au jeu une somme totale $S$, exprimée
en pistoles. Ils jouent à pile
ou face, et le premier qui remporte $n$ parties gagne la somme. Nous
dirions aujourd'hui que Pascal recherche l'espérance de gain $E$ d'un
joueur à un stade quelconque de la partie ou, ce qui revient au
même, sa probabilité $p$ de gagner. Or d'un point de vue technique,
il ne calcule pas une espérance mais un nombre $E'$ qu'il définit comme ``ce
qu'il en appartient sur l'argent de l'autre'', c'est-à-dire comme
l'espérance de gain ôtée de $S/2$ (on a donc $E'=E-S/2$).
Et de même qu'à l'espérance on associe la probabilité par la formule
$E=pS$, Pascal
associe à son nombre $E'$ un coefficient $p'$ tel que
$E'=p'S/2$ ($S/2$ car il ne s'intéresse qu'au parti sur l'argent de l'autre
joueur et non à la somme totale). $p'$ est donc la proportion de
l'argent de l'autre qui revient au joueur ayant l'avantage. Ce n'est pas une 
probabilité de gagner, mais une sorte de coefficient d'avantage sur
son adversaire. 

On peut sans peine passer des concepts aujourd'hui usuels à ceux de
Pascal. De $E'=E-S/2$, on déduit en reportant les valeurs de $E$
et $E'$~: 

\[
p'\frac{S}{2}=pS-\frac{1}{2}S
\]

Ou encore~:

\begin{equation}
p=\frac{p'}{2}+\frac{1}{2}
\end{equation}

Il ne faut pas attacher une trop grande importance à l'usage de $E'$
et $p'$, puisque la notion d'`avantage sur l'adversaire' n'apparaît
nulle part ailleurs que dans la lettre du 29 juillet 1654. Pascal 
semble ne l'avoir retenue que dans la mesure où comme on l'a déjà
vu, elle lui permettait d'écrire une `belle formule', dont nous allons
prouver la véracité.

\subsection*{\citer{44}{3}}

Pascal donne un résultat pour le moins obscur. Il y est question
de ``la valeur de la dernière partie de deux'', de la ``valeur de la
dernière partie de trois'' \ldots 

Le sens du mot `valeur' ne pose pas
problème : à un stade quelconque du jeu, une victoire fait augmenter
l'espérance de gain du joueur, et cette augmentation est la ``valeur'' de
la partie . En revanche, l'expression ``dernière partie de deux''
laisse perplexe. Voici mon interprétation de ce qu'a voulu dire
Pascal.

Au moment où un joueur gagne le jeu, le score est nécessairement de la forme
\score{n}{n-k+1}, avec $k\geq 2$. Pascal dit simplement que la
valeur de la dernière partie jouée est alors de ${S}/{2^{k-1}}$, ce
qui se vérifie par une récurrence sur $k$. à mon avis, quand Pascal dit à Fermat
``dernière partie de 2'', il dit en fait que k vaut 2 ;~quand il parle
de la ``dernière partie de 3'', il considère que $k=3$ etc. Voyons ce
qui me permet d'affirmer cela~:

Admettons que Pascal considère des parties
`canoniques' dans lesquelles, pour arriver au score final de \score{n}{n-k+1}, les
joueurs commencent par faire une quasi-égalité à \score{n-k}{n-k+1}, et
où par la suite le joueur ayant $(n-k)$ remonte son handicap jusqu'à la
victoire finale. On voit alors que c'est bien à une suite de $k$ victoires
consécutives que le vainqueur doit de gagner. Voilà en quoi
l'expression ``dernière partie de 2'' peut correspondre à $k=2$,
``dernière partie de 3'' à $k=3$ et ainsi de suite.

%Evidemment, on peut alors se demander de
%la dernière partie de quoi il est question. Remarquons tout d'abord
%que la première partie a une valeur qui ne dépends que de $n$. Pour la
%dernière, il n'en est rien. En trois parties gagnantes par exemple
%($n=3$), Pascal montre qu'en cas de score 2-0, la valeur d'une
%nouvelle victoire est de 8 pistoles. Mais en cas de score 2-2, cette valeur
%passe à 32 pistoles, puisqu'elle assure le gain à partir d'une seule
%victoire et d'une situation équitable. En fait Pascal part de la fin
%de la partie et considère sans le dire explicitement le nombre de
%victoires du gagnant depuis sa dernière 
%défaite. Voilà pour le sens de cette mystérieuse valeur d'une
%``dernière partie de 3''.

%A ce stade, on peut déjà voir qu'en toute rigueur l'affirmation de
%Pascal est fausse (si j'ai bien compris). En effet, la valeur de la
%dernière partie dépends uniquement de l'espérance de gain du joueur
%avant cette partie et non pas du nombres de victoire consécutives
%qu'il a pu obtenir auparavant.

%Ceci étant dit, l'affirmation de Pascal est (je l'espère) moins
%obscure. Hélas Pascal ne donne aucune indication de
%démonstration. D'après le plan de la lettre, il semble qu'il y soit
%parvenu par des méthodes `conditionelles'. Et en effet, on peut s'y
%l'on en prend la peine suivre ses pas dans le cas général ($n$ joueur,
%somme $S$). 

%Si le score est de $n-1$ à $n-1$, on voit que la valeur sur l'argent
%de l'autre de dernière
%partie est $S/2$. Si le score est $n-1$ à $n-2$, 

\subsection*{\citer{44}{4}\ à \citer{44}{8}}

Pascal donne la valeur de la première partie gagnée dans le cas d'une
partie à n joueurs. Ici le recours à une formulation moderne permet de
gagner en concision~: Pascal veut calculer ``La proportion des
premières parties''. En langage moderne, nous dirions qu'il cherche
l'espérance de gain d'un joueur sachant qu'il a déjà une victoire à
zéro à
son actif.

Pascal dit que (\citer{44}{8})~:

\begin{equation}
    p'=\frac
	{\displaystyle \prod_{i=1}^{n-1} (2i-1)}
	{\displaystyle \prod_{i=1}^{n-1} (2i)}
\end{equation}

\subsection*{\citer{44}{9}\ à \citer{45}{4}}

Pascal `démontre' la propriété ci-dessus. En fait, Pascal ne s'attarde
que sur les parties triviales de la démonstration, laissant dans
l'ombre ce qui est à mon sens le plus important. Au \citerp{44}{9},
Pascal annonce que ses méthodes de probabilités conditionnelles ne lui
ont pas permis de trouver quelque chose d'intéressant dans le cas
général de $n$ joueurs. Les paragraphes 10 à
12 (ce
dernier étant la traduction du latin en français du \citerp{44}{11}) sont entièrement
consacrés à l'énoncé de la formule suivante~:

\begin{equation}
	 	 \comb{2(n-1)}{n-1}
		+\comb{2(n-1)}{n}
		+\comb{2(n-1)}{n+1}
		+\cdots
		+\comb{2(n-1)}{2(n-1)}
	=
		2^{2n-3}+\frac{1}{2}\comb{2(n-1)}{n-1}	
\end{equation}

Pascal ne donne pas la moindre indication sur la manière de prouver
cette formule, qui résulte des Conséquences septième et neuvième du
\emph{Traité du triangle arithmétique}, et il n'en
fait aucun usage explicite par la suite. 

à la  \citer{45}{2}, Pascal remarque qu'après une victoire le jeu se
résout obligatoirement en $2(n-1)$ coups au plus.

Au \citerp{45}{3}, Pascal énonce au moyen du nombre générique le
résultat suivant~: 

\begin{equation}
	p'=\frac{1}{2^{2n-2}}\comb{2n-2}{n-1}
\end{equation}

Ce résultat n'est pas immédiat, et Pascal n'en donne pas la
preuve. J'en propose une qui est probablement proche de
la démonstration originelle puisqu'elle utilise (3)~:

Je mène par \score{1}{0}. Le jeu se décidera donc en au plus $2(n-1)$
coups. Supposons que ces $2(n-1)$ coups soient effectivement joués
\footnote{Cet artifice de calcul, la condition feinte, fait l'objet
d'une
controverse analysée dans la première partie de ce mémoire.}.
Ma probabilité de gagner est égale au 
nombre $N_f$ de cas favorables sur le nombre total $2^{2(n-1)}$ de
cas.
Pour gagner, il suffit que sur les $2^{2(n-1)}$ coups j'en gagne au
moins $n-1$. Donc on peut écrire :

\[
	N_f= 	 \comb{2(n-1)}{n-1}
		+\comb{2(n-1)}{n}
		+\comb{2(n-1)}{n+1}
		+\cdots
		+\comb{2(n-1)}{2(n-1)}
\]

En utilisant (3), on voit donc que :

\[
	N_f=	\frac{1}{2}\comb{2(n-1)}{n-1}+2^{2n-3}
\]

Ainsi :

\[
	p= \frac
		{ \frac{1}{2}\comb{2(n-1)}{n-1}+2^{2n-3} }
		{ 2^{2(n-1)} }
\]

Ou encore :

\[
	p=\frac{1}{2}
			+ \frac{1}{2}
			  \frac
				{\comb{2(n-1)}{n-1}}
				{2^{2(n-1)}}
\]

En utilisant (1) on trouve bien :

\[
	p'=\frac{1}{2^{2n-2}}\comb{2n-2}{n-1}
\]

Il se peut que Pascal ait calculé directement $p'$ sans passer par le
calcul de $p$. Je ne suis pas parvenu à imaginer un tel
raisonnement. Mais dans cette hypothèse,
on pourrait supposer que Pascal a développé au moins provisoirement un
calcul des `coefficients d'avantage sur son adversaire', très commode pour le
calcul pas à pas des espérances de gains (voir la méthode `pas à pas')
mais dont les principes me restent mystérieux s'il fallait l'appliquer
directement à la méthode par les
combinaisons. Cela me paraît toutefois invraisemblable puisque dans
tous ses autres travaux, Pascal n'a plus jamais recours à de telles notions. 

Enfin, au \citerp{45}{4}, Pascal rappelle son résultat (2), sans donner plus
d'information sur sa preuve. Là encore, j'ai quelque idée sur la façon
dont Pascal a pu procéder pour montrer (2) à partir de (4)~:

D'après une expression classique de $C_n^p$, qui apparaît tout à
la fin du \emph{Traité du triangle arithmétique}, on a :

\[
	\frac
		{\comb{2n-2}{n-1}}
		{2^{2n-2}}
	=
	\frac
		{n(n+1)(n+2)\cdots (2n-2)}
		{1(2)(3)\cdots (n-1)}
	\frac{1}{2^{2n-2}}
\]

Ce dont on déduit :

\[
	p'=
	\frac
		{n(n+1)(n+2)\cdots (2n-2)}
		{2(4)(6)\cdots 2(n-1)}
	\frac{1}{2^{n-1}}
\]
	
Donc prouver (2) revient à monter que :

\[
	2^{n-1}(1(3)(5) \cdots (2n-3))=n(n+1)(n+2)\cdots (2n-2)
\]

Cette égalité se montre par récurrence\footnote{à mon avis Pascal
n'a pas démontré cette relation par récurrence. Mais ce point est
secondaire étant donné le caractère élémentaire de la propriété.} On la
vérifie sans peine pour 
des petites valeurs de $n$. Suposons qu'elle soit vraie pour $n$ fixé et
explicitons alors le membre gauche dans le cas $n+1$ :

\[
	2^{n}(1(3)(5) \cdots (2n-1))
	=
	2(2n-1)(n(n+1)(n+2) \cdots (2n-2)) 
\]

On obtient alors ce qu'on voulait montrer :

\[
	2^{n}(1(3)(5) \cdots (2n-1))
	=
	(n+1)(n+2) \cdots (2n-2)(2n-1)(2n).
\]

\newpage

\section{Bibliographie}

\noindent Voici les recueils des textes de Pascal, Fermat et Huygens
	utilisés pour la préparation de ce mémoire. 

\begin{itemize}

\item
\emph{\OE uvres complètes} de Fermat.

On y trouvera une 
correspondance exhaustive avec des indications précises sur la datation
des lettres.

\item
\emph{Précis des \oe uvres mathématiques de
Fermat}, chez Jacques Gabay, 1989.

Curieuse 
compilation (sans table des matières). Sûrement plus facile à se
procurer que les \emph{\OE uvres 
complètes}. En regardant successivement toutes les pages, on trouvera
après un nombre fini d'itérations quelques lettres de Fermat. 

\item
\emph{\OE uvres complètes} de Pascal, au Seuil, 1988.

Comportent tous les traités. On ne trouvera que les lettres de Pascal
lui-même, ce qui n'est pas 
très pratique. Mais les passages en latin sont traduits.

\item {\emph{Pensées}} de Pascal

L'édition du Livre de poche a un index plus copieux que celui des
\emph{\OE uvres complètes}. Permettez-moi d'illustrer ce point par une
anecdote personnelle : un soir,
j'entrepris de lire quelques pensées, Je tombai sur un passage plus
abscons encore que les autres, plus profond peut-être. Au moment
de m'endormir, je réalisai qu'en fait je lisais l'index.

\end{itemize}

\newpage
\section{Actualité de la géométrie du hasard}

\end{document}